\documentclass[pdftex, 10pt, a4paper, oneside]{article}
\usepackage[english]{babel}
\usepackage[utf8]{inputenc} 
\usepackage{indentfirst} 
\usepackage[a4paper,top=4cm,bottom=4cm,left=3cm,right=3cm]{geometry}
\usepackage[numbers]{natbib}
\usepackage{color}

\usepackage{verbatim}
\usepackage{tikz}
\usetikzlibrary{arrows,chains,matrix,positioning,scopes}
\usepackage{amsthm}
\usepackage{amssymb}
\usepackage{amsmath}
\usepackage{mathrsfs}
\usepackage{latexsym}
\usepackage{graphicx}
\DeclareGraphicsExtensions{.pdf,.jpg}
\usepackage{subfigure}
\usepackage{hyperref}
\usepackage{amsfonts}
\usepackage{bm}
\usepackage{algorithm2e}
\usepackage{afterpage}
\usepackage{fancyhdr}
\author{R\'emi Abgrall\footnote{I\lowercase{nstitut f\"ur} M\lowercase{athematik,} W\lowercase{interthurstrasse} 190, CH 8057  Z\lowercase{\"urich}, S\lowercase{witzerland}.} \, and Davide Torlo\footnotemark[1] \footnote{ C\lowercase{orresponding author,  (\href{mailto:davide.torlo@math.uzh.ch}{davide.torlo@math.uzh.ch}).}} }

\def \L {\mathcal{L}}
\def \R {\mathbb{R}}
\def \B {\mathbb{B}}

\theoremstyle{plain}
\newtheorem{teo}{Theorem}[section]

\newtheorem{proposition}[teo]{Proposition}
\theoremstyle{definition}
\theoremstyle{remark}

\newcommand{\disp}{\displaystyle}

\begin{document}
\title{Asymptotic preserving Deferred Correction Residual Distribution schemes}
\date{}
\maketitle

\begin{abstract}
This work aims to extend the residual distribution (RD) framework to stiff relaxation problems. 
The RD is a class of schemes which is used to solve hyperbolic system of partial differential equations. 
Up to our knowledge, it was used only for systems with mild source terms, such as gravitation problems or shallow water equations. 
What we propose is an IMEX (implicit--explicit) version of the residual distribution schemes, that can resolve stiff source terms, without refining the discretization up to the stiffness scale. 
This can be particularly useful in various models, where the stiffness is given by topological or physical quantities, e.g. multiphase flows, kinetic models, viscoelasticity problems.
Moreover, the provided scheme is able to catch different relaxation scales automatically, without losing accuracy. 
The scheme is asymptotic preserving and this guarantees that in the relaxation limit, we recast the expected macroscopic behaviour. 
To get a high order accuracy, we use an IMEX time discretization combined with a Deferred Correction (DeC) procedure, while naturally RD provides high order space discretization.  
Finally, we show some numerical tests in 1D and 2D for stiff systems of equations.
\end{abstract}

{\small

}

\section{Introduction}
In many models, such as kinetic models, multiphase flows, viscoelasticity or relaxing gas flows, we have to deal with hyperbolic systems with relaxation terms. The relaxation term is often led by a parameter $\varepsilon$, the relaxation parameter, that can represent the mean free path, the average distance between two collisions of particles, the time needed to reach the equilibrium between two phases, etc. Expanding these equations asymptotically with respect to $\varepsilon$, one can find the limit equations that describe the average, effective or macroscopic physical behaviour \cite{natalini,Jin_Xin,russo}.

In particular, we focus on the kinetic model proposed by Aregba-Driollet and Natalini in \cite{natalini,zbMATH01519661}. This model is able to solve any hyperbolic system of equation, through an artificial relaxation, which leads to a linear advection system with a relaxation source term. 
It can be used to test classical hyperbolic systems in the relaxation limit case. 
This model must be subjected to a generalization of Whitham's subcharacteristic condition \cite{natalini, Jin_Xin}, which assures that we are adding numerical viscosity to the limit equations. We use this model to approximate transport linear equation, Burgers' equation and Euler equation in 1D and 2D. 
There are various other models and physical problems which behave similarly to this kinetic model. The perspective is, in future, to extend the method to multiphase flows, viscoelasticity problems, and so on.

We use the residual distribution (RD) framework \cite{abgrall, paola_svetlana, mario_enciclopedia, ricc_abgrall_expl_RK} to discretize our space. This class of schemes is a generalization of finite element schemes and allows to recast different well known finite element, finite volume and discontinuous Galerkin schemes \cite{Abgrall2017jcp}. 
The main ingredients of the scheme are three: we have to compute total residuals for each cell of the discretized domain, then, we have to distribute each residual to degrees of freedom of the cell, finally, we sum all contributions at each node. In order to get a high order scheme, the RD is coupled with a Deferred Correction (DeC) iterative method to have computationally explicit schemes \cite{DeC_Abgrall, Dutt_DeC, minion}. It needs two operators: the first one is a low order method, but easy to be inverted, while the second one, must be higher order, but we do not need to solve it directly. The coupling of these two allows to reach the high order through a few iterative intermediate steps. Thanks to this, we can produce a scheme which is fast, high order and stable. Up to our knowledge, RD was utilised only for hyperbolic equations with mild source terms, such as in gravitation problems or shallow water equations, but never on strongly stiff source terms. 

To deal with the stiffness of the relaxation term, we have to introduce some special treatments. An explicit scheme with CFL conditions tuned on the macroscopic regime would, indeed, present instabilities. 
To properly catch the small scale of the microscopic model, one must classically recur to very fine time and space discretizations that are not always feasible in terms of computational time. 
The natural alternative is to use an implicit or semi--implicit formulation, which guarantees the stability of the scheme. We use an IMEX (implicit--explicit) scheme to treat implicitly the relaxation term and explicitly the advection part \cite{Jin_Xin, russo}. 
Nevertheless, we propose a computationally explicit scheme, thanks to some properties of the model. 
Then, we introduce an IMEX discretization for the DeC RD schemes with the details of its implementation. Furthermore, we prove that the new DeC RD IMEX scheme is asymptotic preserving (AP). 
The AP property of a numerical method allows to preserve the asymptotic behaviour of the model from the microscopic to the macroscopic case. 
These schemes solve the microscopic equations, avoiding coupling of different models, and automatically are able to solve the asymptotic macroscopic limit in a robust way. In the appendix, we also provide a proof of the accuracy of the total scheme.

We show the performance of the high order scheme on some tests. In particular, we simulated different examples in 1D and 2D for linear transport equation and Euler equation. Thanks to these results, we validate the accuracy of our method and the capability of shock limiting along discontinuities.

The outline of the manuscript is as follows. In section \ref{sec:kinetic_scheme} we present the kinetic model we want to solve and the conditions under which it is stable.  In section \ref{sec:RD} we describe the RD schemes for the spatial discretization with the DeC high order time discretization. In section \ref{sec:IMEX}, we need to adjust the time discretization according to an IMEX scheme, to deal with stiff source terms and we prove the asymptotic preserving property of the scheme. We show numerical results for 1D and 2D problems in section \ref{sec:test}. Finally, in section \ref{sec:conclusion}, we describe the conclusions and some future investigations that may be done.

\section{Kinetic relaxation model for hyperbolic systems}\label{sec:kinetic_scheme}
In this section, we introduce the kinetic relaxation model presented by D. Aregba-Driollet and R. Natalini in \cite{natalini,zbMATH01519661}. This is a first step to solve general hyperbolic systems of conservation laws via a relaxed system.

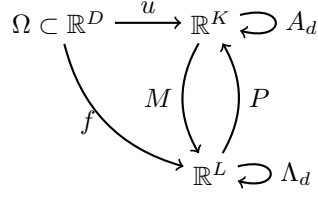
\begin{figure}
	\centering
\begin{tikzpicture}[node distance=2cm] 

  \node (a) {$\Omega\subset\R^D$};
\node[right of=a] (b) {$\R^K$};
\node[below of=b] (c) {$\R^L$};

\draw[thick, ->] (a) -- node[anchor=south]{$u$} (b) ;
\draw[thick, ->] (a) to [bend right] node[anchor=east]{$f$} (c) ;
\draw[thick, ->] (b) to [bend right] node[anchor=east]{$M$} (c) ;
\draw[thick, ->] (c) to [bend right] node[anchor=west]{$P$} (b) ;
\draw[thick, ->] (b) to [loop right] node[anchor=west]{$A_d$} (b) ;
\draw[thick, ->] (c) to [loop right] node[anchor=west]{$\Lambda_d$} (c) ;

\end{tikzpicture}
\caption{Relaxation functions}\label{pic:functions_spaces}
\end{figure}

Let $u: \Omega \subset \R^D \times [0,T] \to \R^K$ be a weak solution of the following system of equations 
\begin{equation}\label{non_linear_natalini}
u_t+	\sum_{d=1}^D \partial_{x_d} A_d(u)=0
\end{equation}
with initial conditions $u(x,0)=u_0(x).$ Here, $A_d:\mathbb{R}^K\to\mathbb{R}^K$ are locally Lipschitz continuous on $\R^K$ with values in $\R^K$. We approximate the problem with a relaxed system
\begin{equation}\label{relaxed_natalini}
f^\varepsilon_t +\sum_{d=1}^D \Lambda_d \partial_{x_d} f^\varepsilon = \frac{1}{\varepsilon} \left( M(Pf^\varepsilon)-f^\varepsilon \right),\qquad \qquad f^\varepsilon(x,0)=f^\varepsilon_0(x)
\end{equation}
where $f^\varepsilon :\Omega \subset \mathbb{R}^D\times [0,T] \to \mathbb{R}^L$ with $M:\mathbb{R}^K\to \mathbb{R}^L$ Lipschitz continuous Maxwellian function, $P:\mathbb{R}^L\to \mathbb{R}^K $ a constant projection matrix ($L>K$) and $\Lambda_d$ diagonal $L\times L $ matrices as sketched in figure \ref{pic:functions_spaces}. 

Moreover, we require that for all $u$ in a certain manifold of interest of $\R^K$ the relations
\begin{equation}\label{kinetic_condition}
\begin{cases}
P(M(u))=u \\
P\Lambda_d M(u)=A_d(u)
\end{cases}
\end{equation}
hold. If $f^\varepsilon$ converges in some strong topology to a limit $f$ and $Pf^\varepsilon_0\to u_0$, then $Pf$ is a solution of the first system \eqref{non_linear_natalini}.
\\
To show this, we define $u^\varepsilon:=Pf^\varepsilon,\, v^\varepsilon_j :=P\Lambda_j f^\varepsilon $ for $j=1,\dots, D$. 
Then we have from \eqref{relaxed_natalini} that
\begin{equation}\label{Pf_equation}
\begin{cases}
\partial_t u^\varepsilon+\sum_{j=1}^D \partial_{x_j} v_j^\varepsilon=0 \vspace{1mm} \\ 
\partial_t v_d^\varepsilon + \sum_{j=1}^D \partial_{x_j} (P\Lambda_j \Lambda_d f^\varepsilon)=\frac{1}{\varepsilon} ( A_d(u^\varepsilon )-v_d^\varepsilon), \quad \forall d \in \lbrace 1, \dots,  D \rbrace
\end{cases}.
\end{equation}
Again, thanks to \eqref{relaxed_natalini}, we consider a formal expansion of $f^\varepsilon$ in Taylor series with respect to $\varepsilon$ in the form of 
\begin{equation}
f^\varepsilon = M(u^\varepsilon) + \varepsilon g^\varepsilon + \mathcal{O}(\varepsilon^2),
\end{equation}
from the second equation of \eqref{Pf_equation} we can write $\forall d = 1,\dots, D$
\begin{align}
v_d^\varepsilon & = A_d(u^\varepsilon ) - \varepsilon \left( 	\partial_t v_d^\varepsilon + \sum_{j=1}^D \partial_{x_j} (P\Lambda_d \Lambda _ j f^\varepsilon)	\right) + \mathcal{O}(\varepsilon^2)\\
&= A_d(u^\varepsilon ) - \varepsilon \left( 	\partial_t v_d^\varepsilon + \sum_{j=1}^D \partial_{x_j} (P\Lambda_d \Lambda _ j M(u^\varepsilon))	\right) + \mathcal{O}(\varepsilon^2).
\end{align}
If we substitute this result in \eqref{Pf_equation}, we get
\begin{align}
\partial_t u^\varepsilon + \sum_{d=1}^D \partial_{x_d} A_d (u^\varepsilon) = \varepsilon \sum_{d=1}^D \partial_{x_d} \left( \partial _t v_d^\varepsilon + \sum_{j=1}^D \partial_{x_j} (P\Lambda_d \Lambda_j M(u^\varepsilon))  \right) + \mathcal{O}(\varepsilon^2).
\end{align}
Now, we have that
\begin{align}
\partial_t v_d^\varepsilon = \partial_t A_d(u^\varepsilon) + \mathcal{O}(\varepsilon) =A_d'(u^\varepsilon)\partial_t u^\varepsilon + \mathcal{O}(\varepsilon) = - \sum_{j=1}^D A'_d(u^\varepsilon) A'_j(u^\varepsilon) \partial_{x_j}u^\varepsilon + \mathcal{O}(\varepsilon) .
\end{align}
Then, we eventually obtain up to second order in $\varepsilon$
\begin{align}
\partial_t u^\varepsilon + \sum_{d=1}^D \partial_{x_d} A_d(u^\varepsilon) = \varepsilon \sum_{d=1}^D \partial_{x_d} \left( \sum_{j=1}^D B_{dj}(u^\varepsilon) \partial_{x_j} u^\varepsilon \right)
\end{align}
where 
\begin{align}
B_{dj}(u):= P\Lambda_d\Lambda_j M'(u) - A'_d(u)A'_j(u)
\end{align}
is a $K\times K $ matrix.\\
This limit equation is stable if the following condition holds:
\begin{align}\label{whit}
\sum_{j,d=1}^D (B_{dj} \xi_j,\xi_d) \geq 0, \qquad \forall \xi_1, \dots, \xi_D \in \R^K .
\end{align}
This property is a generalization of the \textit{Whitham's subcharacteristic condition} \cite{natalini, Jin_Xin,zbMATH01519661}.

We have to choose $M,P,\Lambda $ that respect conditions \eqref{kinetic_condition} to completely define the kinetic model. First of all, let us take in consideration $L=N\times K$ with $P=(I_K,\dots,I_K)$ the juxtaposition of $N$ blocks of identity matrices $I_K\in \mathbb{R}^{K\times K}$. Here, we can consider several $f_n^\varepsilon\in \mathbb{R}^K$ with $n=1,\dots, N$ instead of a single vector $f^\varepsilon\in \mathbb{R}^{N\times K}$, several Maxwellians $M_n:\R^K\to \R^K$ and a block diagonal matrix $\forall d =1,\dots, D$  $$\Lambda_d=diag (C_1^{(d)}, \dots, C_N^{(d)} )\qquad C_n^{(d)}=\lambda_n^{(d)} I_K, \quad \text{with }\lambda_n^{(d)} \in \mathbb{R}, \, \forall n=1,\dots, N.$$ With this formalism we can rewrite \eqref{relaxed_natalini} as
\begin{equation}\label{DRM}
\begin{cases}
\partial_t f_n^\varepsilon + \sum_{d=1}^D \lambda_n^{(d)} \partial_{x_d} f_n^\varepsilon = \frac{1}{\varepsilon} \left(  M_n(u^\varepsilon) -f_n^\varepsilon\right),\qquad \forall n =1,\dots , N \\
u^\varepsilon=\sum_{n=1}^N f_n^\varepsilon
\end{cases}.
\end{equation}

Let us present the \textit{diagonal relaxation method (DRM)}. Here $N=D+1$. Then we have to define Maxwellians $M_n$ and matrices $C_j^{(d)}$. Take $\lambda>0$, that will be chosen according to \textit{Whitham's subcharacteristic condition }\eqref{whit}, and 
\begin{equation}
C_j^{(d)}=\begin{cases} -\lambda I_K  &j=d
\\
\lambda I_K  &j=D+1\\
0  &\text{else}
\end{cases}.
\end{equation}
The Maxwellians can be defined as follows:
\begin{equation}
\begin{cases}
M_{D+1}(u)=\left( u + \frac{1}{\lambda	} \sum_{d=1}^D A_d(u) \right) /(D+1)\\
M_j(u)=-\frac{1}{\lambda} A_j(u)+M_{D+1}(u)
\end{cases}
\end{equation}

For one--dimensional system of conservation laws this formulation coincides with Jin--Xin relaxation model \cite{Jin_Xin}, the simplest example that we can think of in this context. Indeed, if we set $u^\varepsilon :=Pf^\varepsilon$ and $v^\varepsilon:=P\Lambda f^\varepsilon$, we get
\begin{equation}\label{DRM_Jin_Xin}
\begin{cases}
\partial_t u^\varepsilon +\partial_x v^\varepsilon =0\\
\partial_t v^\varepsilon +\partial_x u^\varepsilon = \frac{1}{\varepsilon} (A(u^\varepsilon)-v^\varepsilon).
\end{cases}
\end{equation}

\section{Residual distribution schemes}\label{sec:RD}
Let us now introduce the spatial and time discretization given by RD schemes \cite{ abgrall99,  mario_enciclopedia} and DeC approach \cite{DeC_Abgrall, Dutt_DeC}. 

\subsection{Notation}
Let us start introducing the notation of RD schemes.
For sake of simplicity, we explain the RD approach for steady equations, the time derivative part will be discussed in section \ref{sec:Dec}. So, we can focus on the following equation $$\nabla \cdot A(U) -S(U)=0 .$$
We define a triangulation $\Omega_h$ on our domain $\Omega$ and denote by $K$ the generic element of the mesh and by $h$ the characteristic mesh size (implicitly supposing some regularity on the mesh). 
\begin{figure}[ht!]
\begin{center}
\includegraphics[scale=0.56]{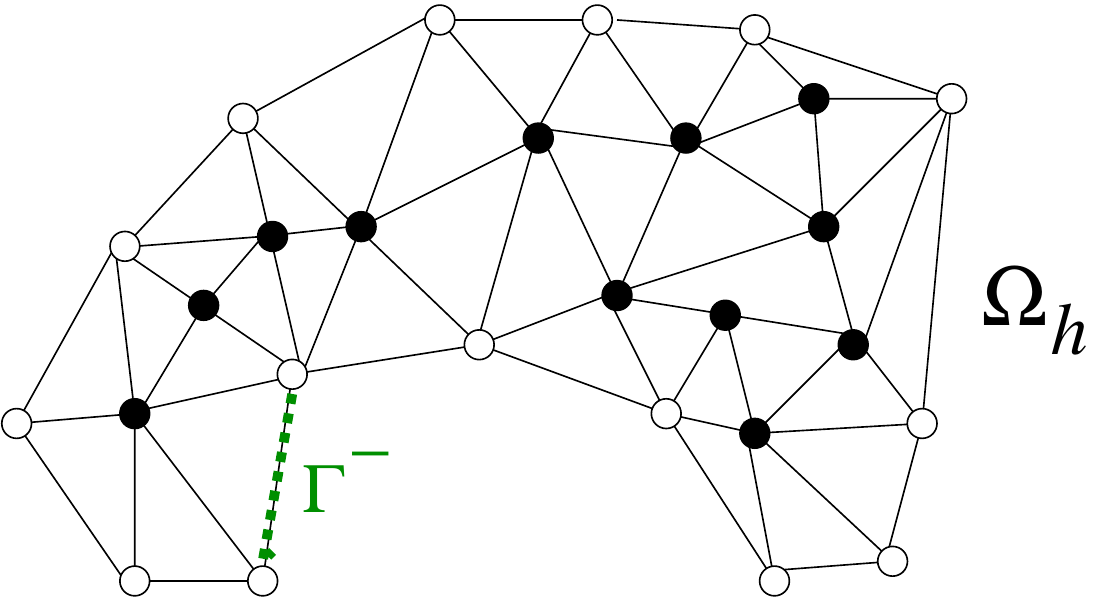}
\end{center}\caption{Triangulation of the domain $\Omega$}\label{triangulation}
\end{figure}
Following the ideas of the Galerkin finite element method (FEM), we use a solution approximation space $V_h$ given by globally continuous piecewise polynomials of degree $k$:
\begin{equation}\label{discretized_functional_space}
V_h=\lbrace U \in \mathcal{C}^0(\Omega_h), \, U|_K\in \mathbb{P}^d,\,  \forall K \in \Omega_h \rbrace.
\end{equation}
Now we can rewrite the numerical solution $U_h(x) \approx U(x) $ as a linear combination of basis functions $\varphi_\sigma \in V_h$:
\begin{equation}\label{approximation_function}
U_h(x)=\sum\limits_{\sigma \in D_h} U_\sigma \varphi_\sigma(x) = \sum_{K\in \Omega_h} \sum_{\sigma \in K} U_\sigma \varphi_\sigma |_K(x), \quad \forall x \in \Omega
\end{equation}
where $D_h$ is the set of all the degrees of freedom of $\Omega_h$, so that $ \lbrace \varphi_\sigma : \sigma \in \mathcal{D}_h \rbrace $ is a basis for $V_h$, and the coefficient $U_\sigma$ must be found by a numerical method.\\
\subsection{Residual distribution scheme}
\begin{figure}
\begin{center}
\includegraphics[width=0.32\textwidth, trim={0 60  0 0}, clip]{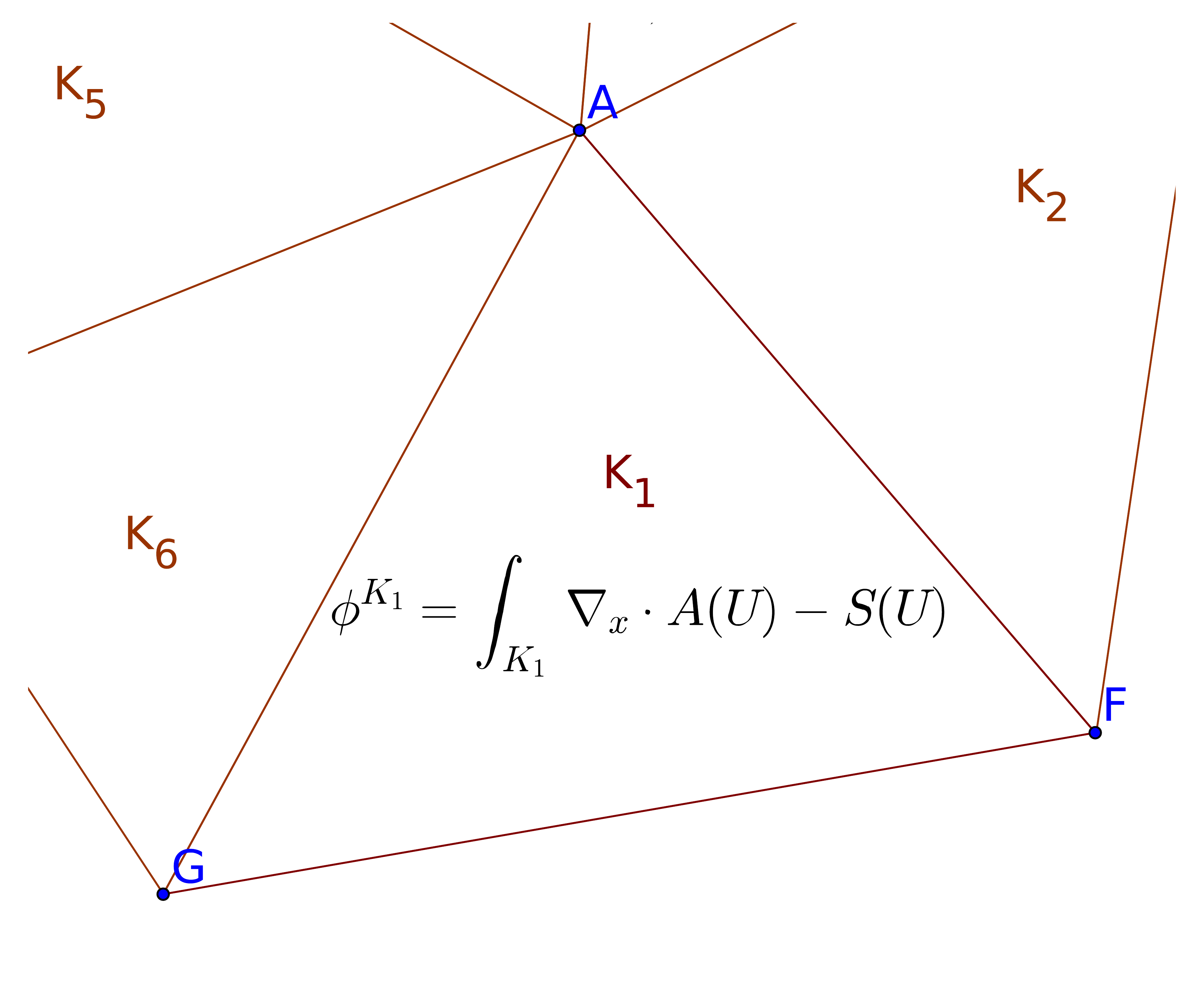}  \hfill
\includegraphics[width=0.32\textwidth, trim={0 60  0 0}, clip]{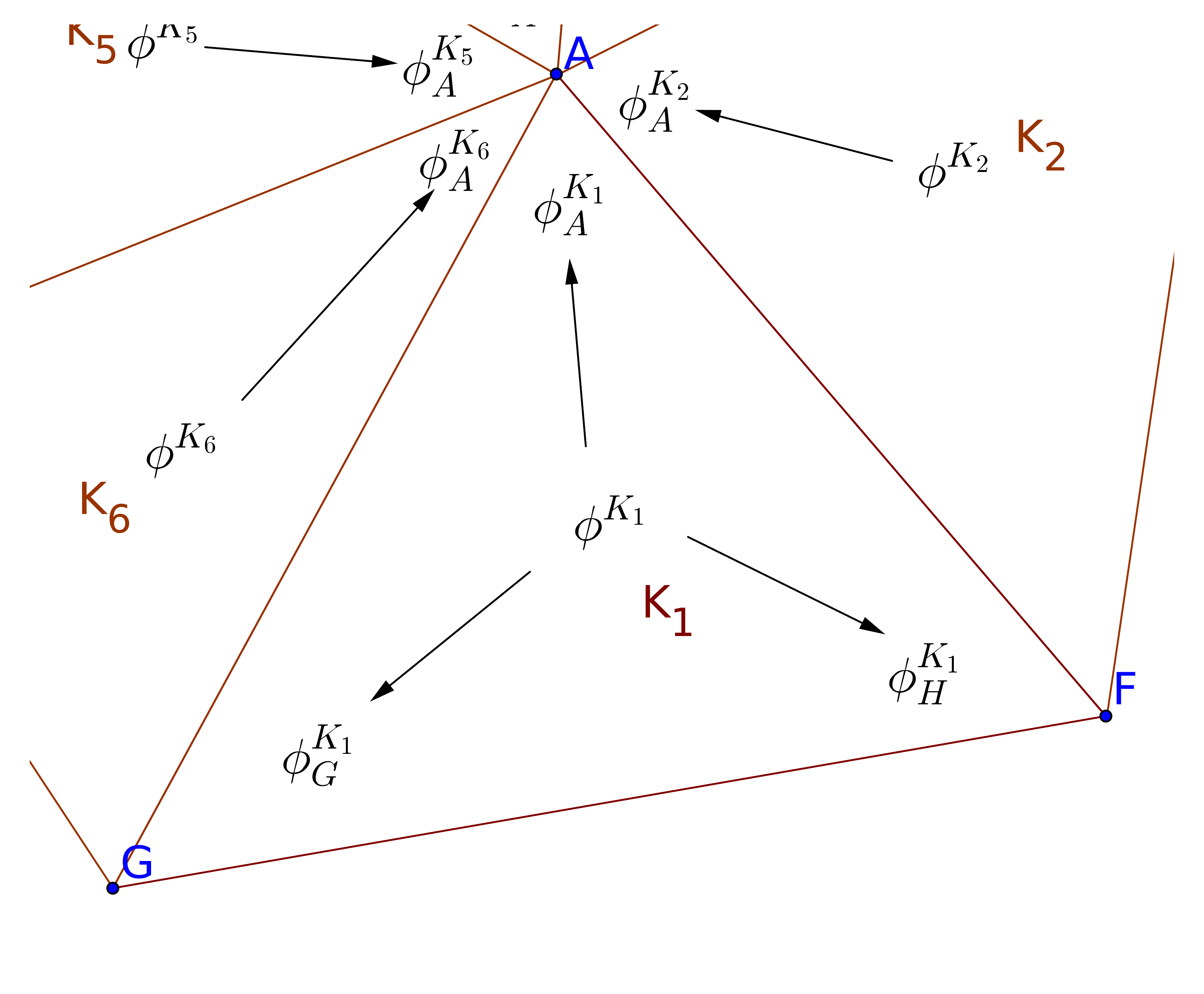} \hfill 
\includegraphics[width=0.32\textwidth, trim={0 0  0 0}, clip]{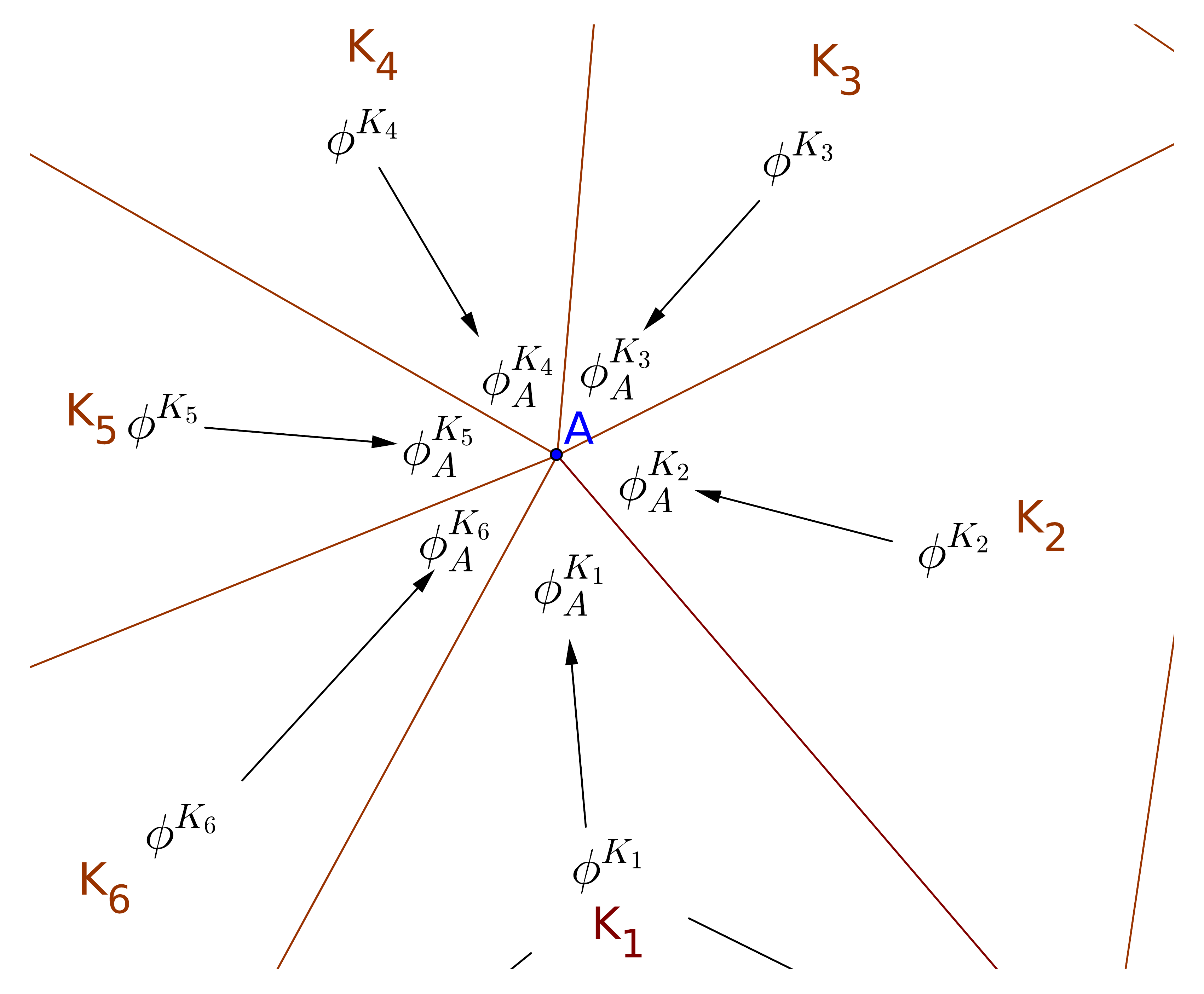} 
\end{center}\caption{Defining total residual, nodal residuals and building the RD scheme}\label{picture_RD}
\end{figure}

RD schemes can be summarized as follows.
\begin{enumerate}
\item Define $\forall K \in \Omega_h$ a fluctuation term (total residual) \begin{equation}
\phi^K=\int_K \left ( \nabla \cdot A(U_h) - S(U_h) \right )dx
\end{equation}
\item Define a nodal residual $\phi^K_\sigma$ as a contribution to fluctuation term $\phi^K$ for each degree of freedom $\sigma$ within the element $K$, so that the sum of all the contributions over an element is the fluctuation itself, i.e., 
\begin{equation}
\phi^K=\sum_{\sigma\in K} \phi^K_\sigma, \quad \forall K \in \Omega_h. 
\end{equation}
In appendix \ref{app:residual_distribution} or \cite{zbMATH05020825,zbMATH05920225} one can find more details on possible definitions of the nodal residuals.

\item The resulting scheme is obtained by summing all the nodal residual contributions of one degree of freedom from different elements $K$, that is 
\begin{equation}
\sum_{K|\sigma \in K } \phi^K_\sigma =0, \quad \forall \sigma \in D_h.
\end{equation}
This is a RD scheme.
\end{enumerate}
The main sketch of the scheme is done in picture \ref{picture_RD}.
The key of the scheme is the definition of nodal residuals. This choice is leading the whole spatial discretization. The scheme can be highly accurate in space, just choosing higher order polynomial basis functions and consistent nodal residuals. In \cite{abgrall99, DeC_Abgrall, Abgrall2017jcp} it has been shown that well known finite element or finite volume schemes (such as SUPG, DG, FV-WENO, etc.) can be rewritten in terms of RD, just choosing the proper nodal residuals.

Details and some examples of the schemes can be found in the appendix \ref{app:residual_distribution}.

\subsection{Time discretization}\label{sec:Dec}
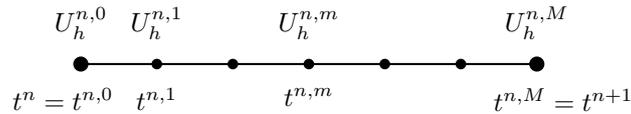
\begin{figure}[h]
	\centering
\begin{tikzpicture}
\draw [thick]   (0,0) -- (6,0) node [right=2mm]{};
\fill[black]    (0,0) circle (1mm) node[below=2mm] {$t^n=t^{n,0} \,\, \quad$} node[above=2mm] {$U_h^{n,0}$}
                            (1,0) circle (0.7mm) node[below=2mm] {$t^{n,1}$} node[above=2mm] {$U_h^{n,1}$}
                            (2,0) circle (0.7mm) node[below=2mm] {}
                            (3,0) circle (0.7mm) node[below=2mm] {$t^{n,m}$} node[above=2mm] {$U_h^{n,m}$}
                            (4,0) circle (0.7mm) node[below=2mm] {}
                            (5,0) circle (0.7mm) node[below=2mm] {}
                            (6,0) circle (1mm) node[below=2mm] {$\qquad t^{n,M}=t^{n+1}$} node[above=2mm] {$U_h^{n,M}$};
\end{tikzpicture}
\caption{Subtimesteps}\label{subtimestep}
\end{figure}

For time discretization, we want to get a high order accurate approximation. To do so, we discretize the timestep $[ t^n, t^{n+1} ]$ into $M$ subtimesteps $[t^{n,0},t^{n,1}],\dots,[t^{n,M-1}, t^{n,M}]$ and the variable $U_h$ in time at each subtimestep $U^{n,m}_h$ as in picture \ref{subtimestep}.

Using the Picard--Lindel\"of theorem, we can write for $m=1,\dots ,M$
\begin{equation}\label{eq:picard}
U_h^{n,m} - U_h^n+\int\limits_{t^n}^{t^{n,m}} \left( \nabla \cdot A(U_h(x,s))-S(U_h(x,s)) \right )ds=0.
\end{equation}
For sake of simplicity, we drop from now on the subscript $h$. More precisely, the scheme that we want to solve is a system of equations, where each entry is the discretization of \eqref{eq:picard} for a different $m=1,\dots, M$. In practice, we can write it as
{\small
\begin{equation}\label{eq:L2}
\begin{split}
\L^2_\sigma &(U^{n,0},\dots,U^{n,M})=\\
=&\begin{pmatrix}
\displaystyle\sum\limits_{K|\sigma\in K} \displaystyle\int_K \varphi_\sigma (U^{n,1}-U^{n,0}) dx + \sum\limits_{K|\sigma\in K} \displaystyle \int\limits_{t^{n,0}}^{t^{n,1}}\mathcal{I}_M (\phi^K_\sigma (U^{n,0}),\dots, \phi^K_\sigma (U^{n,M}),s)  ds    \\
\vdots\\
\displaystyle\sum\limits_{K|\sigma\in K}  \displaystyle\int_K \varphi_\sigma (U^{n,M}-U^{n,0}) dx +  \sum\limits_{K|\sigma\in K} \disp \int\limits_{t^{n,0}}^{t^{n,M}}\mathcal{I}_M (\phi^K_\sigma (U^{n,0}),\dots, \phi^K_\sigma (U^{n,M}),s)  ds   
\end{pmatrix}.
\end{split}
\end{equation}  }
Here, we have $M$ equations with $M$ unknowns $U^{n,1},\dots, U^{n,M}$, $\mathcal{I}_M$ is an interpolation polynomial in nodes $\lbrace t^{n,m}\rbrace_{m=0}^M $ and the time integration is computed using quadrature formulas in the interpolation points. Of course, this system may contain a lot of nonlinear terms as functions of $U$, so we would like not to solve it directly. Nevertheless, the solution to \eqref{eq:L2} is what we are interested in. It is an approximation of the real solution with an accuracy of order $M+1$ in time and $d+1$ in space, where $d$ is the degree of utilised polynomials.

The spirit of the DeC algorithm is to use two schemes, one high order and another one explicit or easy to solve.
So, we introduce a first order approximation of the scheme $\L^2$, that we will call $\L^1$:
{\small
\begin{equation}\label{eq:L1}
\begin{split}
\L^1_\sigma& (U^{n,0},\dots,U^{n,M})=\\
=&\begin{pmatrix}
\displaystyle (U^{n,1}_\sigma-U^{n,0}_\sigma)\sum\limits_{K|\sigma\in K}  \int_K \varphi_\sigma  dx + \sum\limits_{K|\sigma\in K}   \displaystyle\int\limits_{t^{n,0}}^{t^{n,1}}\mathcal{I}_0 (\phi^K_\sigma (U^{n,0}),\dots, \phi^K_\sigma (U^{n,M}),s)   ds    \\
\vdots\\
\displaystyle (U^{n,M}_\sigma-U^{n,0}_\sigma)\sum\limits_{K|\sigma\in K} \int_K \varphi_\sigma  dx + \sum\limits_{K|\sigma\in K} \displaystyle\int\limits_{t^{n,0}}^{t^{n,M}}\mathcal{I}_0 (\phi^K_\sigma (U^{n,0}),\dots, \phi^K_\sigma (U^{n,M}),s)   ds   
\end{pmatrix}.
\end{split}
\end{equation}
}

The first simplification we applied is a mass lumping on the derivative in time, substituting $U$ with $U_\sigma$. This is only possible if $|\mathcal{C}_\sigma|=\sum_K \int_K \varphi_\sigma(x) dx > 0.$ For this reason, we will always consider Bernstein polynomials $\B^d$, which are nonnegative everywhere, instead of Lagrange polynomial $\mathbb{P}^d$.

The second one is in the residual part, where we substituted the high order interpolant $\mathcal{I}_M$ with a piecewise constant interpolant $\mathcal{I}_0$, which is explicit or easy to solve. An example of interpolant polynomial can be  $\mathcal{I}_0(\phi^K_\sigma (U^{n,0}),\dots, \phi^K_\sigma (U^{n,M}),s) \equiv \phi^K_\sigma (U^{n,0})$. The detail of the interpolant will be given in section \ref{sec:IMEX}. The approximation error brought from these two approximations is a $\mathcal{O}(\Delta t + \Delta x)$.

\subsection{Deferred Correction algorithm}
Now, we present the deferred correction (DeC) algorithm to couple the two formulations. It was introduced by A. Dutt in \cite{Dutt_DeC} and we can see another approach in \cite{minion}, but we follow the formulation by Abgrall in \cite{DeC_Abgrall}. The aim of DeC schemes is to avoid implicit methods, without losing the high order of accuracy of a  scheme. In our case, the high order method that we want to approximate is $\L^2$ of \eqref{eq:L2}. To use the DeC procedure, we also need another method, which is easy and fast to be solved, we use diagonal mass matrix explicit methods, with low order of accuracy $\L^1$, as in \eqref{eq:L1}. The DeC algorithm is providing an iterative procedure that wants to approximate the solution of the $\L^2$ scheme $U^*$ in the following way.
\begin{equation}\label{DeC_method}
\begin{split}
&\L^1(U^{(1)})=0, \\
&\L^1(U^{(k)})=\L^1(U^{(k-1)})-\L^2(U^{(k-1)}) \text{ with }k=2,\dots,K,
\end{split}
\end{equation}
where $K$ is the number of iterations that we compute. In particular, we need as many iteration as the order of accuracy that we want to reach: $K=d+1=M+1$. Notice that, in every step, we solve the equations for the unknown variable $U^{(k)}$ which appears only in the $\L^1$ formulation, the one that can be solved easily. While $\L^2$ is only applied to already computed predictions of the solution $U^{(k-1)}$. Thus, we can state the following proposition as  in \cite{DeC_Abgrall}.
\begin{proposition}\label{DeC_prop}
Let $\L^1$ and $\L^2$ be two operators defined on $\mathbb{R}^m$, which depend on the discretization scale $\Delta\sim \Delta x \sim \Delta t$, such that
\begin{itemize}
\item $\L^1$ is coercive with respect to a norm, i.e., $\exists \alpha_1 >0$ independent of $\Delta$, such that for any $U,V$ we have that $$\alpha_1||U-V||\leq ||\L^1 (U)-\L^1 (V)||,$$
\item $\L^1 - \L^2$ is Lipschitz with constant $\alpha_2>0$ uniformly with respect to $\Delta$, i.e., for any $U,V$
$$
||(\L^1(U)-\L^2(U))-(\L^1(V)-\L^2(V))||\leq \alpha_2 \Delta ||U-V||.
$$
\end{itemize}
We also assume that there exists a unique $U^*_\Delta$ such that $\L^2(U^*_\Delta)=0$. Then, if $\eta:=\frac{\alpha_2}{\alpha_1}\Delta<1$, the DeC is converging to $U^*$ and after $k$ iterations the error $||U^{(k)}-U^*||$ is smaller than $\eta^k||U^{(0)}-U^*||$.

\end{proposition}
The proof of the proposition can be found in appendix \ref{DeC_proposition}, while the proof of the properties of $\L^1$ and $\L^2$, which depend on their definitions, can be found for our specific case in appendix \ref{Dec_properties}.

The theorem tells us that, if the method $\L^2$ is accurate with order of accuracy $r$, then we should perform $r$ iterations for every timestep of the method and that we need only $r - 1$ sub-time steps. 
For example, if we use $\mathbb{B}^1$ basis functions, we will have 2 iterations of the DeC method (1 prediction and 1 correction) with 1 sub-time steps ($t^{n,0}=t^n,\,t^{n,1}=t^{n+1}$): this amounts to one version of the second order Runge Kutta method, see \cite{ricc_abgrall_expl_RK}. 
For $\mathbb{B}^2$, we need 3 iterations (1 prediction, 2 corrections) and 2 sub-time steps ($t^{n,0}=t^n,\,t^{n,1}=\frac{1}{2}(t^n+t^{n+1}),\,t^{n,2}=t^{n+1}$) and so on. If not specified, in all our test cases we will use the same number of degree of polynomial, corrections-1 and subtimesteps, i.e., $d=K-1=M$.

\section{IMEX asymptotic preserving kinetic scheme}\label{sec:IMEX}
Before introducing an IMEX scheme, let us explain what is the problem concerning the kinetic model that we are considering. Solving equation \eqref{relaxed_natalini}, we have to be careful in treating the source term. If we discretize it in an explicit way, it would produce strongly stiff terms as $\varepsilon \to 0$. 
To classically solve this problem, one should take very small $\Delta t $ values of the order of $\Delta t \sim \varepsilon$. 
On the other hand, the solution of the system would induce very long computational time. 
That is why, this can not always be a feasible way. The alternative is to treat implicitly the source term. Namely, we can use this type of time discretization:
\begin{align}\label{IMEX_time_discretization}
&\frac{f^{n+1,\varepsilon}-f^{n,\varepsilon}}{\Delta t} +\sum_{d=1}^D \Lambda_d \partial_{x_d} f^{n,\varepsilon} = \frac{1}{\varepsilon} \left( M(Pf^{n+1,\varepsilon})-f^{n+1,\varepsilon} \right), \\  &f^{0,\varepsilon}(x)=f^\varepsilon_0(x),
\end{align}
where the superscript index in $f^n$ indicates the $n$-th timestep. This type of discretization is called implicit--explicit (IMEX), since the advection term is explicit, while the source term is implicit. This approach guarantees stability to the time discretization and we can relax the constraint on $\Delta t$ until the usual CFL conditions proportional to the eigenvalue of the jacobian of the flux, which is $\lambda$ in DRM model. Overall, the time-step can be chosen such that $\Delta t \leq \text{CFL } \lambda \Delta x $, where the CFL depends on the degree of the used polynomial basis functions.

As it is written, the time discretization \eqref{IMEX_time_discretization} presents some nonlinear implicit terms. We can get rid of this technical problem, so that the scheme turns out to be computationally explicit. What we can notice is that the source is depending nonlinearly on $Pf^{n+1,\varepsilon} = u^{n+1,\varepsilon} $ and linearly on $f^{n+1,\varepsilon}$. To reach our goal, we can solve the following auxiliary equation for $u^{n+1,\varepsilon}$, which is the results of the multiplication of \eqref{IMEX_time_discretization} by $P$ and properties  \eqref{kinetic_condition}:
\begin{equation}\label{method_limit_update}
\frac{u^{n+1,\varepsilon} -u^{n,\varepsilon }}{\Delta t}+ \sum_{d=1}^D P \Lambda _d \partial_{x_d} f^{n,\varepsilon}=0.
\end{equation} 
We can see that for this equation we are simply applying forward Euler method, which is explicit, since the source term turns out to be zero.
So, we can solve it and then substitute $u^{n+1,\varepsilon}$ in equation \eqref{IMEX_time_discretization} and solve it without recurring to implicit methods nor inversion of mass matrices. 
Indeed, the equation \eqref{IMEX_time_discretization} can be rewritten in the following form, where the right--hand--side is explicit:
\begin{equation}\label{IMEX_time_discretization_2}
f^{n+1,\varepsilon} \left( \frac{1}{\Delta t} + \frac{1}{\varepsilon} \right) =\frac{f^{n,\varepsilon}}{\Delta t} -\sum_{d=1}^D \Lambda_d \partial_{x_d} f^{n,\varepsilon} + \frac{1}{\varepsilon} M(u^{n+1,\varepsilon}) .
\end{equation}
One can, indeed, express the variable $f^{n+1}$ in the following way
\begin{equation}\label{IMEX_time_discretization_3}
f^{n+1,\varepsilon}  = \frac{\varepsilon}{\Delta t+ \varepsilon}f^{n,\varepsilon} -\frac{\varepsilon \Delta t}{\Delta t+ \varepsilon}\sum_{d=1}^D \Lambda_d \partial_{x_d} f^{n,\varepsilon} + \frac{\Delta t}{\Delta t +\varepsilon} M(u^{n+1,\varepsilon}).
\end{equation}
We can see that, in this formulation, $\varepsilon$ does not appear alone in any denominator, so, for $\varepsilon\to 0$, $f^{n+1,\varepsilon}$ is well defined and tends to the Maxwellian $M(u^{n+1, \varepsilon})$.

\subsection{Residual distribution IMEX operators} \label{sec:DeC_IMEX}
What we need to do now, is to apply the IMEX time discretization to the DeC and RD frameworks. This implies the change of the time discretization only of the operator $\L^1$. Indeed, that is the only operator that we actually need to invert to get solutions of the DeC algorithm. While, we can not modify $\L^2$ because we do not want to drop the order of accuracy and because it will be anyway computed on previously computed solutions.  

To do so, we want to choose the zero order interpolant $\mathcal{I}_0$ in a way that the source term is evaluated constantly on the end of the subtimestep, namely in $t^{n,m}$, while the advection term is evaluated on the beginning of the timestep $t^{n,0}$, i.e., 
\begin{equation}
\begin{split}
\mathcal{I}_0(\phi^K_\sigma (f^{n,0}),\dots, \phi^K_\sigma (f^{n,M}),s) \equiv \phi_{ad,\sigma}^K (f^{n,0}) + \phi_{source, \sigma}^K (f^{n,m}).
\end{split}
\end{equation}
This requires a further definition of the nodal residuals that splits the source term and the advection part. The choice of the source residual is done accordingly to IMEX discretization. Indeed, what we require is its implicitness, the linear dependence on $f^{n,m}_\sigma$ and that it does not depend on other degrees of freedom. To reach these goals, we will perform a mass lumping on the whole source term and we evaluate everything in $t^{n,m}$. This results in 
\begin{equation}
\begin{split}
\phi^K_{source, \sigma}= \int_K \varphi_\sigma(x) \frac{M(Pf^{n,m,\varepsilon}_\sigma)-f^{n,m,\varepsilon}_\sigma}{\varepsilon} dx.
\end{split}
\end{equation}
This allows us to collect $f^{n,m,\varepsilon}_\sigma$ on the left hand side of the equation and solve it explicitely. The advection part $\phi_{ad,\sigma}^K$ can be defined in different ways \cite{mario_enciclopedia, paola_svetlana, abgrall99}. We give some examples in appendix \ref{app:residual_distribution}.  Anyway, in this time discretization, it will be always explicit.

From now on we will drop the index $n$ that indicates the timestep we are referring to and the index $\varepsilon$ which refers to relaxation variables. They will be used only when necessary.

Overall, if we define $|\mathcal{C}_\sigma|:=\int_\Omega \varphi_\sigma(x) dx$, the $\L^1$ operator will be at the $m$--th component
\begin{subequations}\label{eq:L1_IMEX_operator}
\begin{equation}\label{eq:L1_u_operator}
\L^{1,m}_{\sigma,u} (f^{0}, u^{m})=  |\mathcal{C}_\sigma| (u^{m}_\sigma  - Pf^{0}_\sigma) +\Delta t^m \sum_{K|\sigma \in K}  P \phi_{ad,\sigma}^K(f^{0});
\end{equation}
\begin{equation}\label{eq:L1_f_operator}
\begin{split}
\L^{1,m}_\sigma (f^{0}, f^{m})= & |\mathcal{C}_\sigma| \left( 1+\frac{\Delta t^m}{\varepsilon} \right) f^{m}_\sigma  - |\mathcal{C}_\sigma|f^{0}_\sigma +\\
&+ \Delta t^m \sum_{K|\sigma \in K} \phi_{ad,\sigma}^K(f^{0}) - |\mathcal{C}_\sigma|\frac{\Delta t^m}{\varepsilon} M(u^{m}_\sigma).
\end{split}
\end{equation}
We can see that both the equations of the $\L^1$ with the IMEX discretization are computationally explicit. Moreover, as before, we can see that, as $\varepsilon\to 0$, equation \eqref{eq:L1_f_operator} does not lead to terms with $\varepsilon$ alone at the denominator. Indeed, it can be rewritten as
\begin{equation}\label{eq:L1_f_operator_explicit}
\begin{split}
\L^{1,m}_\sigma (f^{0}, f^{m})= & f^{m}_\sigma  - \frac{\varepsilon}{\varepsilon + \Delta t^m} f^{0}_\sigma +\\
&+ \frac{\varepsilon \Delta t^m}{|\mathcal{C}_\sigma| (\varepsilon + \Delta t^m)}\sum_{K|\sigma \in K} \phi_{ad,\sigma}^K(f^{0}) - \frac{\Delta t ^m}{\varepsilon + \Delta t^m} M(u^{m}_\sigma).
\end{split}
\end{equation}
\end{subequations}
Finally, we can write a general term of the correction DeC procedure for the $(k+1)$th correction and the $m$th subtimestep. First, we have the $u$ auxiliary equation
\begin{subequations}\label{eq:DeC_IMEX_operator}
\begin{equation}\label{eq:DeC_u_operator}
\begin{split}
\L^{1,m,(k+1)}_{\sigma,u}&  -\L^{1,m,(k)}_{\sigma,u} +\L^{2,m,(k)}_{\sigma,u} =|\mathcal{C}_\sigma| (u^{m,(k+1)}_\sigma  - u^{m,(k)}_\sigma) +\\
  + \sum_{K|\sigma \in K} \bigg[& \int_K \varphi_\sigma(x) (u^{m,(k)}(x)-u^{0,(k)}(x) )dx +\\ +&\int_{t^{n,0}}^{t^{n,m}} \mathcal{I}_M(P \phi_{\sigma,ad}^K(f^{0,(k)}),\dots, P \phi_{\sigma,ad}^K(f^{M,(k)}) ,s)ds \bigg];
  \end{split}
\end{equation}
and, then, the $f$ equation
{\small
\begin{equation}\label{eq:DeC_f_operator}
\begin{split}
&\L^{1,m,(k+1)}_{\sigma} -\L^{1,m,(k)}_{\sigma} +\L^{2,m,(k)}_{\sigma} = \\
& |\mathcal{C}_\sigma| \left( 1+\frac{\Delta t^m}{\varepsilon} \right) (f^{m,(k+1)}_\sigma -  f^{m,(k)}_\sigma )  - |\mathcal{C}_\sigma|\frac{\Delta t^m}{\varepsilon} \left(M\left(u^{m,(k+1)}_\sigma \right) -M\left(u^{m,(k)}_\sigma\right)\right) +\\
&+ \sum_{K|\sigma \in K} \left[ \int_K \varphi_\sigma(x) (f^{m,(k)}(x)-f^{0,(k)}(x) )dx +\int_{t^{0}}^{t^{m}} \mathcal{I}_M( \phi_{\sigma}^K(f^{0,(k)}),\dots,  \phi_{\sigma}^K(f^{M,(k)}) ,s)ds \right] .
\end{split}
\end{equation}}
\end{subequations}
Again, thanks to the factor $\left( 1+\frac{\Delta t^m}{\varepsilon} \right)$ in front of the unknown, we are sure not to have any stiff term, even in the source of $\L^2$.

\subsection{AP property}\label{sec:AP_proof}
An \textit{asymptotic preserving (AP)} scheme preserves the asymptotic behaviour of the model from the microscopic to the macroscopic case. It solves the microscopic equations, avoiding coupling of different models, and, automatically, it is able to solve the asymptotic macroscopic limit as the relaxation parameter $\varepsilon$ tends to its limit.

\begin{figure}[ht!]
	\centering
	
	\begin{tikzpicture}[node distance=2.5cm] 
	
	\node (a) {$\mathcal{F}^\varepsilon_\Delta$};
	\node[right of=a] (b) {$\mathcal{F}^0_\Delta$};
	\node[below of=b] (c) {$\mathcal{F}^0$};
	\node[below of=a] (d) {$\mathcal{F}^\varepsilon$};
	
	\draw[thick, ->] (a) -- node[anchor=south]{$\varepsilon\to 0$} (b) ;
	\draw[thick, ->] (d) -- node[anchor=south]{$\varepsilon\to 0$} (c) ;
	\draw[thick, ->] (b) -- node[anchor=west]{$\Delta \to 0$} (c) ;
	\draw[thick, ->] (a) -- node[anchor=east]{$\Delta \to 0$} (d) ;
	\end{tikzpicture}	
	
	\caption{Asymptotic preserving schemes}\label{ap_figure}
\end{figure}
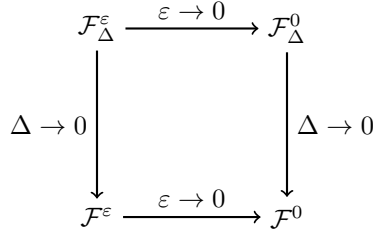

The behaviour of an AP scheme is sketched in figure \ref{ap_figure}. Let us call $\mathcal{F}^\varepsilon$ the microscopic model which depends on $\varepsilon$ and its asymptotic macroscopic limit $\mathcal{F}^0:=\lim\limits_{\varepsilon\to 0} \mathcal{F}^\varepsilon$. 
We denote the numerical discretization of $\mathcal{F}^\varepsilon$ as $\mathcal{F}_\Delta^\varepsilon$, where $\Delta$ is the mesh size and/or the time step length (in our case they are always linked by some CFL conditions). 
Then, we call the asymptotic limit as $\varepsilon\to 0 $ of this scheme $\mathcal{F}^0_\Delta:=\lim\limits_{\varepsilon\to 0} \mathcal{F}_\Delta^\varepsilon$ (for fixed $\Delta$), if it exists. We can say that the scheme $\mathcal{F}^\varepsilon_\Delta$ is an AP scheme, if $\mathcal{F}^0_\Delta$ is a consistent and stable approximation of $\mathcal{F}^0$, i.e.,  $\mathcal{F}_\Delta^0=\mathcal{F}^0+\mathcal{O}(\Delta)$. 
In our model, the limit model $\mathcal{F}^0$ is the equation \eqref{non_linear_natalini}  and the relaxed model $\mathcal{F}^\varepsilon$ is the equation \eqref{relaxed_natalini}. 
What we need to check is that the discrete model, namely, the scheme that we proposed, is asymptotic preserving. 
This implies that, first, we let $\varepsilon\to 0$, then the discretization scale $\Delta\to 0$. 
In other words, we can consider $\frac{\varepsilon}{\Delta}=o(1)$.

To start, let us suppose that the initial conditions $u^\varepsilon_0$ and $f^\varepsilon_0$ verify $f^\varepsilon_0=M(u_0^\varepsilon)$. We can prove that this is also true for the beginning of every timestep by induction. We want to show that at the end of each timestep we maintain the following relation for $u^\varepsilon=Pf^\varepsilon$ at the discrete level:
\begin{align}
&\frac{ u^{n+1} - u^{n}}{\Delta t} +\sum_{d=1}^D\partial_{x_d} A_d(u^{n+1}) + \mathcal{O}(\varepsilon)+ \mathcal{O}(\Delta ) =0. \label{eq:AP_1st}
\end{align}
To prove it by induction, we want to add the following relation in the induction hypothesis
\begin{align}
&\frac{ f^{n+1} - f^{n}}{\Delta t} +\sum_{d=1}^D\partial_{x_d} \Lambda_d(f^{n+1}) -\frac{M(u^{n+1})-f^{n+1}}{\varepsilon}+ \mathcal{O}\left(\frac{\Delta }{\varepsilon}\right)+ \mathcal{O}(\Delta ) =0, \label{eq:AP_3rd}
\end{align}
that implies
\begin{align}
&f^{n+1} = M(u^{n+1}) + \mathcal{O}(\varepsilon)+\mathcal{O}(\Delta ). \label{eq:AP_2nd}
\end{align}
The initial conditions verify the hypothesis \eqref{eq:AP_2nd} for $n=0$. So we check the $n+1$th timestep given the relations \eqref{eq:AP_1st} and \eqref{eq:AP_2nd} for the $n$th timestep. We start from the prediction $\L^1=0$, in forms \eqref{eq:L1_u_operator} and \eqref{eq:L1_f_operator}. Since the scheme begins at each step with \eqref{eq:L1_u_operator}, we can write $\forall m \in [1,\dots,M]$: 
\begin{subequations}
\begin{align}
&\frac{u^{m,(1)}_\sigma-u^{0}_\sigma}{\Delta t^m} +\frac{1}{|\mathcal{C}_\sigma |}\sum_{K|\sigma\in K}P\phi^K_{ad,\sigma}(f^{0})=0,\label{eq:AP_u_discretization0}
\end{align}
and, if we use the fact that the sum of nodal residual is a consistent discretization of space derivatives as shown in \cite{DeC_Abgrall}, we get
\begin{align}
&\frac{u^{m,(1)}_\sigma-u^{0}_\sigma}{\Delta t^m} + \sum_{d=1}^D \partial_{x_d} P \Lambda_d  f^{0}_\sigma + \mathcal{O}(\Delta ) =0. \label{eq:AP_u_discretization} 
\end{align}
Using the induction hypothesis on property \eqref{eq:AP_2nd}, we obtain
\begin{align}
&\frac{u^{m,(1)}_\sigma-u^{0,\varepsilon}_\sigma}{\Delta t^m} + \sum_{d=1}^D \partial_{x_d}  P \Lambda_d M(u^{0}_\sigma) + \mathcal{O}(\Delta ) + \mathcal{O}(\varepsilon) =0, \label{eq:AP_u_properties1}
\end{align}
while, using the properties in \eqref{kinetic_condition}, the equation \eqref{eq:AP_u_discretization0} itself and the fact that $A_d$ are Lipschitz continuous, we reach
\begin{align}
&\frac{u^{m,(1)}_\sigma-u^{0}_\sigma}{\Delta t^m}+ \sum_{d=1}^D \partial_{x_d} A_d(u^{m,(1)}_\sigma) + \mathcal{O}(\Delta ) + \mathcal{O}(\varepsilon) =0.\label{eq:AP_u_properties2}
\end{align}
\end{subequations}
 Then, from \eqref{eq:L1_f_operator} we can recast the second property \eqref{eq:AP_3rd}, with similar reasoning:
\begin{subequations}
\begin{align}
0=&\frac{f^{m,(1)}_\sigma -f^{0}_\sigma}{\Delta t^m} + \frac{1}{|\mathcal{C}_\sigma|}\sum_{K|\sigma \in K} \phi^K_{ad,\sigma} (f^{0}) - \frac{ M(u^{m,(1)}_\sigma)- f^{m,(1)}}{\varepsilon},\label{eq:AP_f_discretization}
\end{align}
then, using consistency of residuals, we can say that
\begin{align}
0= & \frac{f^{m,(1)}_\sigma -f^{0}_\sigma}{\Delta t^m} + \sum_{d=1}^D \partial_{x_d} \Lambda_d f^{0}_\sigma - \frac{ M(u^{m,(1)}_\sigma)- f^{m,(1)}_\sigma}{\varepsilon} + \mathcal{O} (\Delta), \label{eq:AP_f_1}
\end{align}
and, finally, substituting \eqref{eq:AP_f_discretization} in \eqref{eq:AP_f_1}, we get 
\begin{align}
0= & \frac{f^{m,(1)}_\sigma -f^{0}_\sigma}{\Delta t^m} + \sum_{d=1}^D \partial_{x_d} \Lambda_d f^{m,(1)}_\sigma - \frac{ M(u^{m,(1)}_\sigma)- f^{m,(1)}_\sigma}{\varepsilon} + \mathcal{O} (\Delta)+ \mathcal{O} \left(\frac{\Delta}{\varepsilon} \right).\label{eq:AP_f_2}
\end{align}
\end{subequations}
We proved that the prediction is asymptotic preserving, since it recast the limit equation \eqref{non_linear_natalini}. A more rigorous proof of a similar property for the norm convergence of the kinetic scheme is in \cite{natalini,zbMATH01519661}. 

What is left to prove are the same properties for every correction $(k+1)$, using induction hypothesis on the previous correction $(k)$. For prediction $(k)=(1)$, we have already given the proof. Now, let us consider $u$ equation in \eqref{eq:DeC_u_operator} to prove property \eqref{eq:AP_1st}.
\begin{subequations}
\begin{equation}\label{eq:AP_correction_u}
\begin{split}
&\L^{1,m}_{\sigma,u} -\L^{1,m}_{\sigma,u} +\L^{2,m}_{\sigma,u} \\
=& \frac{u^{m,(k+1)}_\sigma  - u^{m,(k)}_\sigma}{\Delta t^m} + \sum_{K|\sigma \in K}  \int_K \varphi_\sigma(x) \frac{u^{m,(k)}(x)-u^{0,(k)}(x)}{|\mathcal{C}_\sigma|\Delta t^m} dx +\\
+&\sum_{K|\sigma \in K} \frac{1}{|\mathcal{C}_\sigma|\Delta t^m}\int_{t^{0}}^{t^{m}} \mathcal{I}_M(P \phi_{\sigma,ad}^K(f^{0,(k)}),\dots, P \phi_{\sigma,ad}^K(f^{M,(k)}) ,s)ds=0,
\end{split}
\end{equation}
Then, we apply a mass lumping of time derivative term in $\L^2$, moreover, we know that the quadrature of the interpolant is a first order approximation of any of its points, which are a consistent approximation of the flux. So,
\begin{align}\label{eq:AP_correction_u_1}
\begin{split}
&\frac{u^{m,(k+1)}_\sigma  - u^{m,(k)}_\sigma}{\Delta t^m} +   \frac{u^{m,(k)}_\sigma-u^{0,(k)}_\sigma}{\Delta t^m} +\mathcal{O}(\Delta) +\\
+&\sum_{d=1}^D\partial_{x_d} P\Lambda_d f^{m,(k)}_\sigma+\mathcal{O}(\Delta)=0, \end{split}
\end{align}
we can now apply property \eqref{eq:AP_2nd} in the induction hypothesis on correction $(k)$ and properties \eqref{kinetic_condition} to get 
\begin{align}
\label{eq:AP_correction_u_3}
 \frac{u^{m,(k+1)}_\sigma  - u^{0,(k)}_\sigma}{\Delta t^m} &+ \sum_{d=1}^D \partial_{x_d} A_d(u^{m,(k)}_\sigma) +\mathcal{O}(\Delta)+\mathcal{O}(\varepsilon)=0,
 \end{align}
 and then we can substitute \eqref{eq:AP_correction_u} in \eqref{eq:AP_correction_u_3} to gain another $\mathcal{O}(\Delta)$ using also the Lipschitz continuity of fluxes $A_d$: 
 \begin{align}
 \label{eq:AP_correction_u_4}
 \frac{u^{m,(k+1)}_\sigma  - u^{0,(k+1)}_\sigma}{\Delta t^m} &+ \sum_{d=1}^D \partial_{x_d} A_d(u^{m,(k+1)}_\sigma) +\mathcal{O}(\Delta)+\mathcal{O}(\varepsilon)=0.
\end{align}
\end{subequations}

Then, to prove property \eqref{eq:AP_2nd}, we can proceed from \eqref{eq:DeC_f_operator}. We can split the three terms of the sum $\L^{1,m,(k+1)}_{\sigma} -\L^{1,m,(k)}_{\sigma} +\L^{2,m,(k)}_{\sigma} $. Let us start from $\L^{2,m,(k)}_{\sigma}$:
\begin{subequations}
\begin{align}
\begin{split}\label{eq:AP_correction_L2_k}
\L^{2,m,(k)}_{\sigma} &=   \frac{1}{|\mathcal{C_\sigma}|}\sum_{K|\sigma \in K}\int_K \varphi_\sigma(x) \frac{f^{m,(k)}(x)-f^{0,(k)}(x)}{\Delta t ^m}dx +\\
&+ \frac{1}{|\mathcal{C_\sigma}|\Delta t ^m} \sum_{K|\sigma \in K}\int_{t^{0}}^{t^{m}} \mathcal{I}_M( \phi_{\sigma}^K(f^{0,(k)}),\dots,  \phi_{\sigma}^K(f^{M,(k)}) ,s)ds  
\end{split}
\end{align}
Then, we use a mass lumping on the time derivative, which brings an error of the order of $\Delta$, the fact that the interpolant is a first order approximation of any of the interpolation points and that the residuals are consistent approximation of the flux and the source. So, we obtain
\begin{align}
\L^{2,m,(k)}_{\sigma}&=\frac{f^{m,(k)}_\sigma-f^{0,(k)}_\sigma}{\Delta t ^m} + \sum_{d=1}^D \partial_{x_d} \Lambda_d f^{m,(k)}_\sigma+ \frac{M(u^{m,(k)}_\sigma) - f^{m,(k)}_\sigma}{\varepsilon} + \mathcal{O}(\Delta). \label{eq:AP_correction_L2_k1}
\end{align}
If we then use the induction hypothesis on $(k)$ correction, we get
\begin{align}
\L^{2,m,(k)}_{\sigma}&=\mathcal{O}\left(\frac{\Delta}{\varepsilon}\right) + \mathcal{O}(\Delta). \label{eq:AP_correction_L2_k2}
\end{align}
\end{subequations}

Analogously, for $\L^{1,m,(k)}_\sigma$ we can prove that it is an $ \mathcal{O}\left(\frac{\Delta}{\varepsilon}\right) + \mathcal{O}(\Delta)$ using the induction hypothesis. Finally, using what we just proved, we have that
\begin{subequations}
{\small
\begin{align}\label{eq:AP_correction_f_1}
&\L^{1,m,(k+1)}_{\sigma} -\L^{1,m,(k)}_{\sigma} +\L^{2,m,(k)}_{\sigma} =\L^{1,m,(k+1)}_{\sigma} +\mathcal{O}\left(\frac{\Delta}{\varepsilon}\right) + \mathcal{O}(\Delta)=0.
\end{align}}
If we express explicitly the formula, we get
{\begin{align}\begin{split}
&   \frac{f^{m,(k+1)}_\sigma-f^{0,(k+1)}_\sigma}{\Delta t ^m} + \sum_{K|\sigma \in K} \frac{\phi_{ad,\sigma}^K(f^{0,(k+1)})}{|\mathcal{C_\sigma}|} \\
+&\frac{M(u^{m,(k+1)}_\sigma) - f^{m,(k+1)}_\sigma }{\varepsilon}+\mathcal{O}\bigg(\frac{\Delta}{\varepsilon}\bigg) + \mathcal{O}(\Delta)=0, \end{split} \label{eq:AP_correction_f_2}
\end{align}}
Using the fact that the residuals are a consistent approximation of the fluxes and that the term at the $m$th subtimestep is an approximation of the term at the $0$th time step, up to an  $\mathcal{O}\left(\frac{\Delta}{\varepsilon}\right) + \mathcal{O}(\Delta)$ from \eqref{eq:AP_correction_f_2}, we finally reach
{\small \begin{align}
&  \frac{f^{m,(k+1)}_\sigma-f^{0,(k+1)}_\sigma}{\Delta t ^m} + \sum_{d=1}^D \partial_{x_d} \Lambda_d f^{m,(k+1)}_\sigma +\frac{M(u^{m,(k+1)}_\sigma) - f^{m,(k+1)}_\sigma }{\varepsilon}+\mathcal{O}\left(\frac{\Delta}{\varepsilon}\right) + \mathcal{O}(\Delta)=0 .\label{eq:AP_correction_f_3}
\end{align}}
\end{subequations}
So, we proved property \eqref{eq:AP_3rd} for all subtimesteps and corrections. This implies that the scheme is AP as $\varepsilon\to 0$ for any discretization scale $\Delta$.

\section{Numerical simulations}\label{sec:test}
To validate the scheme we presented, we test the method on different problems. We will show 1D and 2D test cases for scalar equations and systems of equations. Generally, we will start from the asymptotic limit $u$ and we will draw from that the whole kinetic systems for variable $f$. The shown results are related to the variable $u$. We have some parameter to choose in order to perform our tests. First of all, the convection coefficient $\lambda$, which should satisfy the Whitham's subcharacteristic conditions \eqref{whit} and the relaxation parameter $\varepsilon$ that will be often very small to get the asymptotic behaviour. Then, the CFL conditions, namely a bound on the size of $\Delta t$. Thanks to the scheme presented, we do not need CFL conditions linked to the source term, so, we can just choose them such that 
\begin{equation}
\Delta t \leq \frac{\text{CFL } \Delta x}{\lambda},
\end{equation}
 where $\lambda$ is the convection parameter. While, with a standard RD DeC method without IMEX technique, the $\Delta t$ should scale as $$\Delta t \leq \text{CFL} \min \left\lbrace  \frac{\Delta x}{\lambda}, \frac{\varepsilon }{\lambda} \right\rbrace, $$ which would require very small timesteps that lead to a huge computational demand. The CFL number depends on the degree of the polynomial chosen, and it scales as $\frac{1}{d}$, but for a comparison of the methods, we will choose it uniformly through different polynomial degrees. In all our computations we will also specify the $\theta_k$ parameter, which are leading the stabilization of the jump of the derivative, that we are using in the definition of the nodal residual. More details about the used nodal residual and the jump stabilization can be found in appendix \ref{app:residual_distribution} and in \cite{paola_svetlana}.
\subsection{1D numerical tests}
\subsubsection{Burgers' equation} 
First of all, we start with 1D scalar equations. We want to approximate the Burgers' equation, i.e.,
\begin{equation}
\begin{split}
\partial_t u(x,t)+\partial_x \left(\frac{u(x,t)^2}{2}\right)=0,\quad x \in [ 0,1 ], \quad  t \in [0,T]
\end{split}
\end{equation}
using the relaxation system \eqref{relaxed_natalini}. As initial condition, we take $u_0(x)=\sin(2\pi x)$ and $f_0(x)=M(u_0(x))$ and the boundary conditions are periodic. To satisfy Whitham's condition, we choose $\lambda=2$, so that $|A'(u)|=|u|\leq \lambda$ in an area of interest.\\
In following figures some approximated solutions for different number of elements are shown. To solve the equation we used the scheme \eqref{eq:scheme5} in appendix \ref{app:residual_distribution} with $\theta_1=1$ and, only for $\B^3$, we used $\theta_2=0.5$. The relaxation parameter is set to $\varepsilon=10^{-9}$ and CFL $=0.1$. Final time is $T=0.5$.
\begin{figure}[ht!]
\begin{center}
\subfigure[$N=32$]{\includegraphics[width=0.49\textwidth, trim={30 30 30 20},clip] {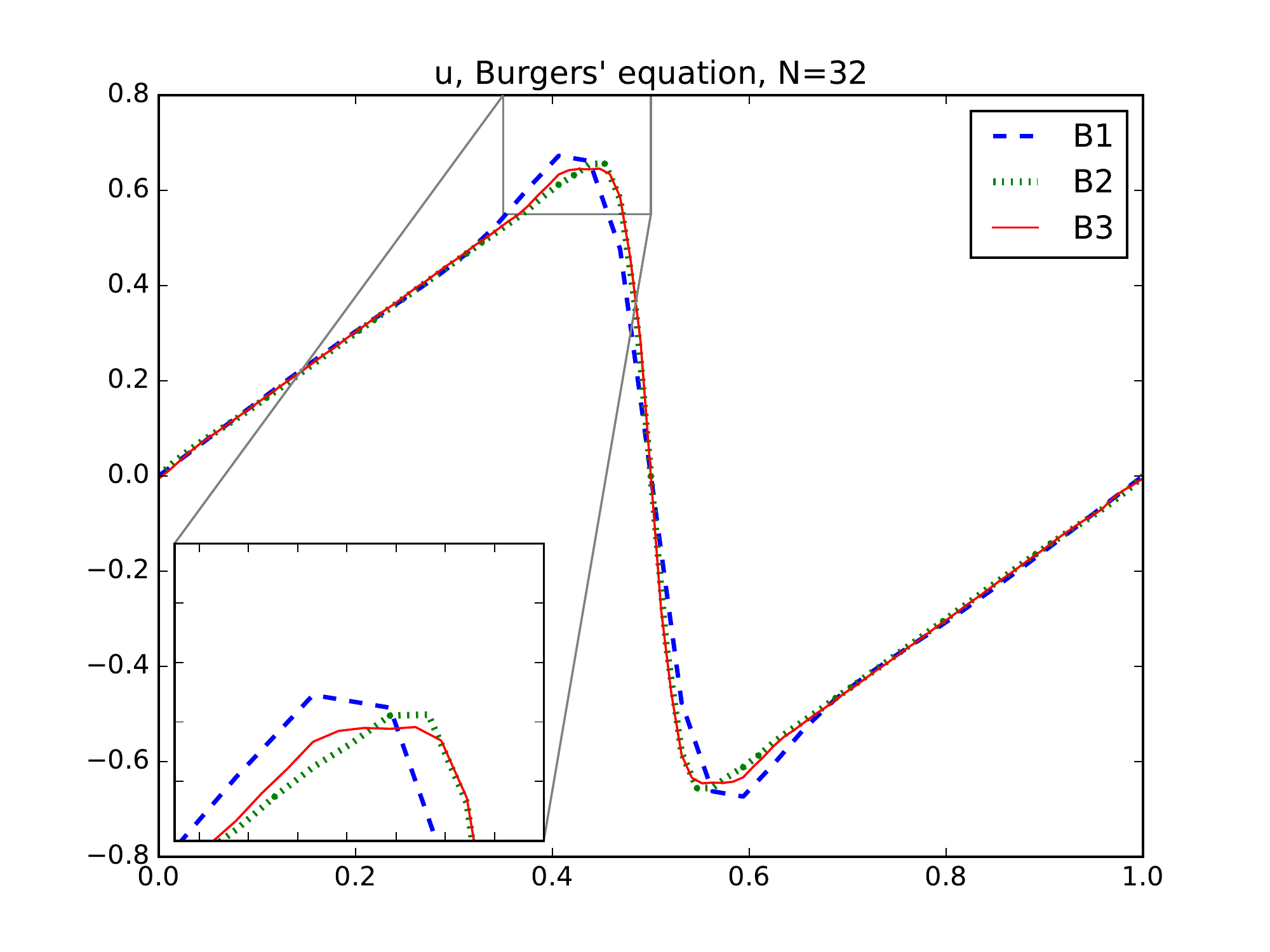} }
\subfigure[$N=128$]{\includegraphics[width=0.49\textwidth, trim={30 30 30 20},clip] {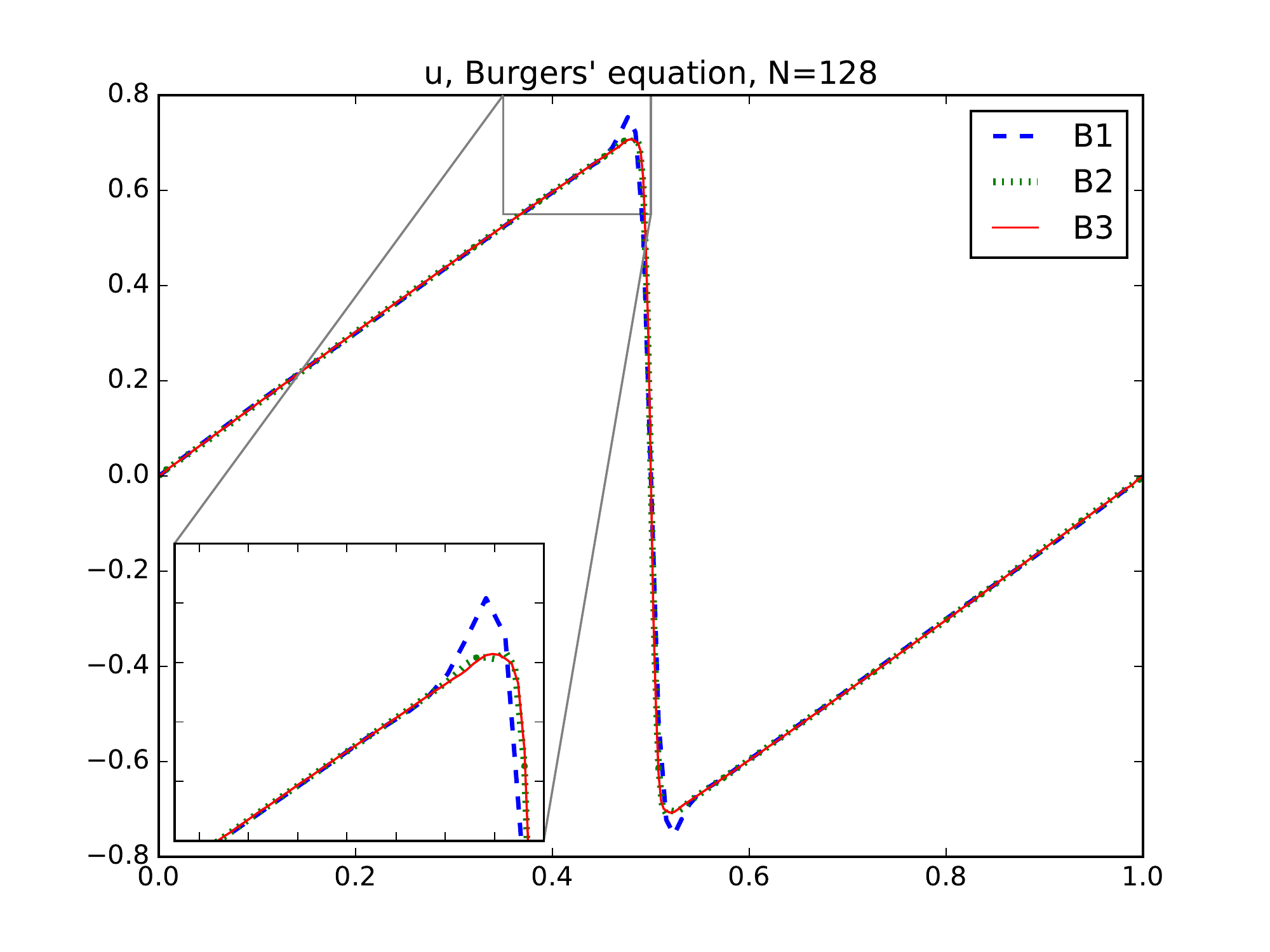}}
\end{center}
\caption{Burgers' equations}\label{picture:burger}
\end{figure}
We can see in picture \ref{picture:burger} that the scheme is well catching the shock position and, as the order of the polynomials increases, we can see improvements in the sharpness of the solution.
\subsubsection{Convergence for linear transport equation}
Then, we test our scheme with different orders to check the convergence rate. 
For all the smooth test cases, where we want to study the order of convergence, we use the scheme which involves only Galerkin residuals and stabilizations of jumps in derivative, as presented in \cite{burman} and in the scheme \eqref{eq:scheme4} in appendix \ref{app:residual_distribution}. We use a linear scalar transport equation $u_t+u_x=0$ as limit equation with the relaxation system presented above,  on domain $[ 0,1 ]$. The initial condition is $u_0(x)=e^{-80(x-0.4)^2}$ and $f_0=M(u_0)$, until final time $T=0.12$ with periodic boundary conditions. 

\begin{figure}[ht!]
\begin{center}
	\subfigure[Scalar 1D convergence]{\includegraphics[width=0.48\textwidth,trim={15 10 15 20},clip]{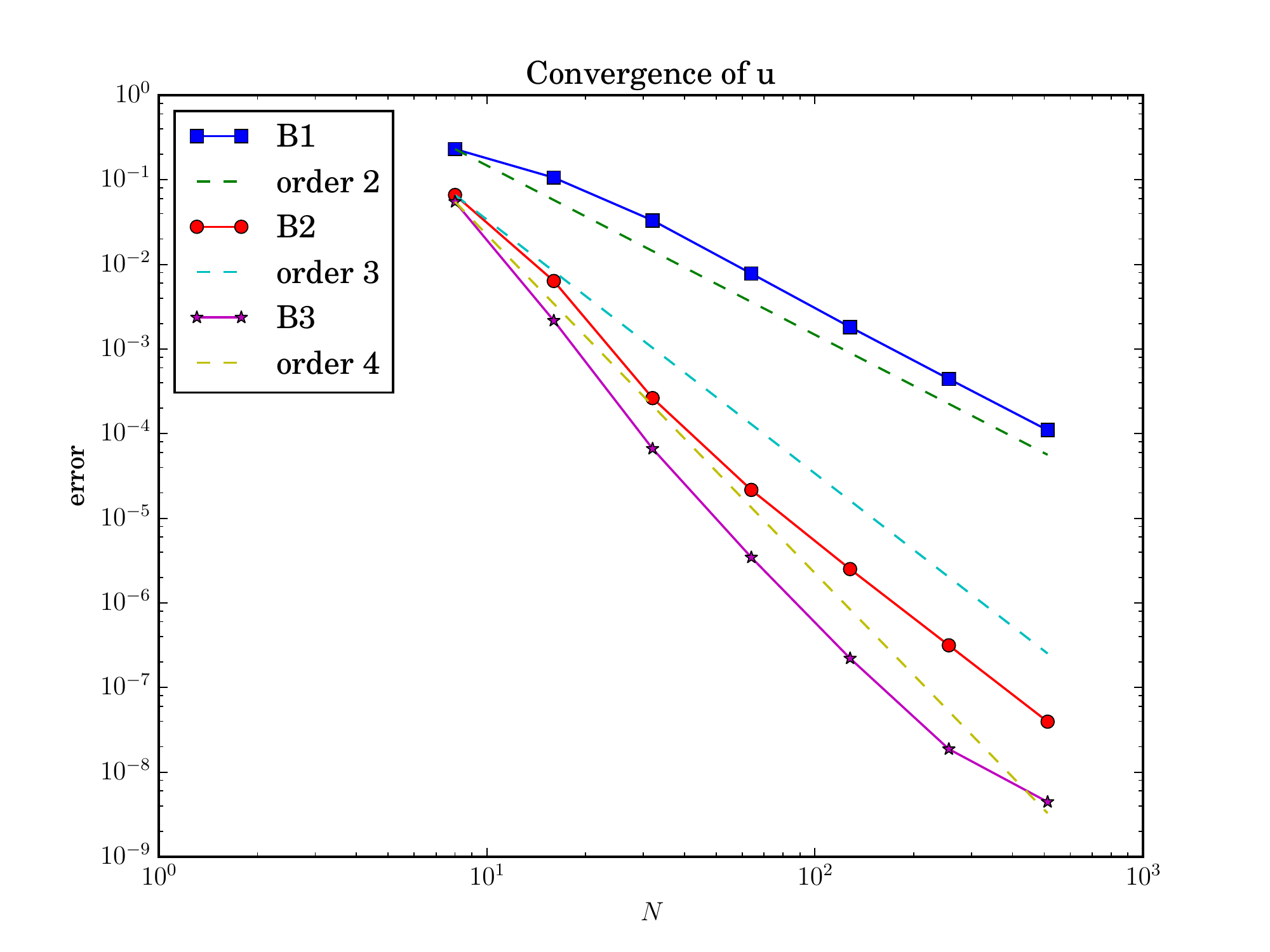} \label{pic:1Dscalar_convergence}}
\subfigure[Varying relaxation parameter]{\includegraphics[width=0.48\textwidth,trim={15 10 15 20},clip]{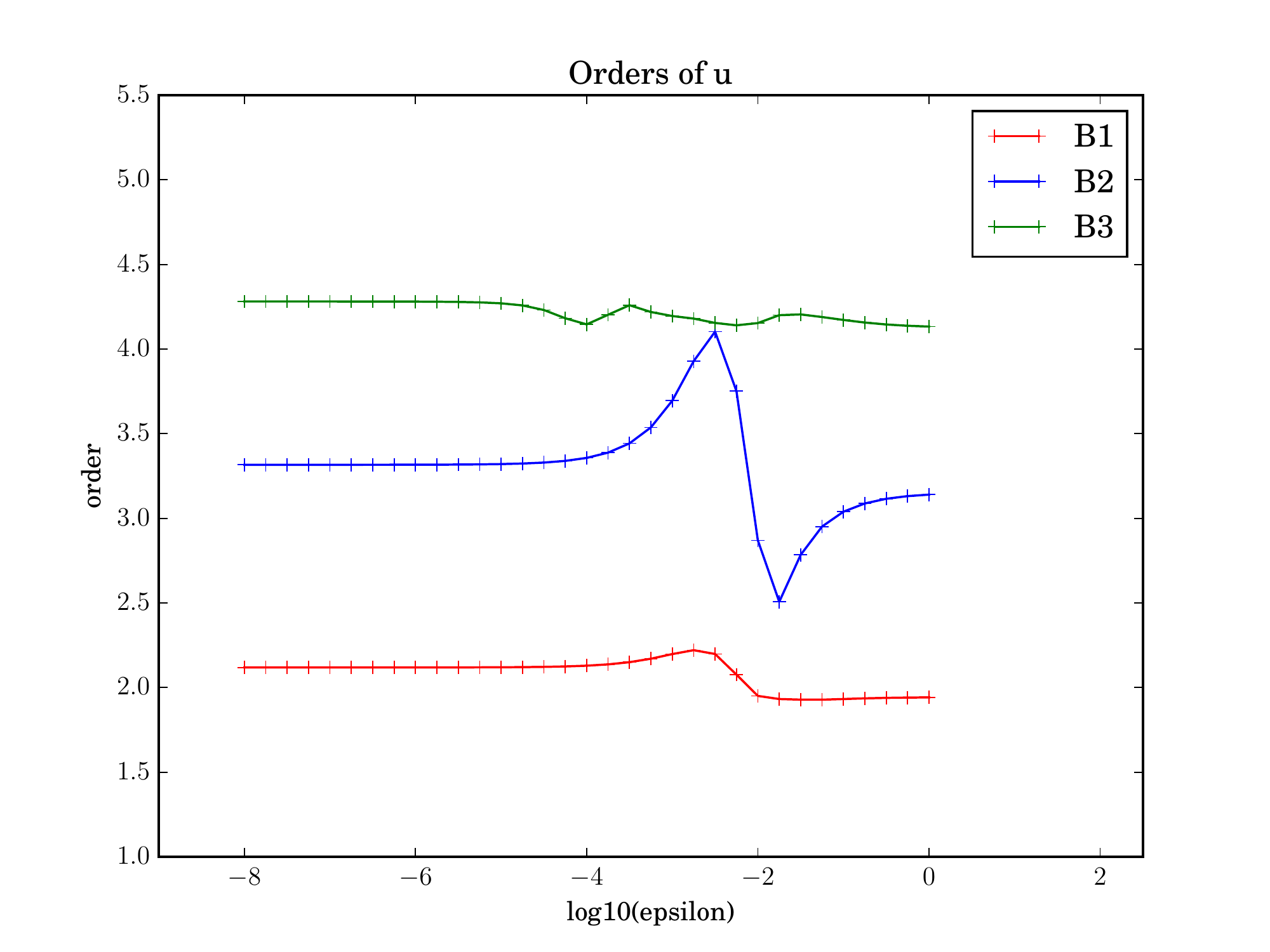} \label{pic:eps_vary} } \caption{Scalar linear 1D test}
\end{center}
\end{figure}

We use the relaxation coefficient $\varepsilon=10^{-9}$, convection $\lambda=1.5$ and CFL=$0.1$. In particular, for $\B^1$ we used $\theta_1 = 1$, for $\B^2$ we used $\theta_1=1,\, \theta_2=0$ and for $\B^3$ we used $\theta_1=1,\,\theta_2=5$. Final time of the solution is $T=0.12$. For $\B^3$ we see that only increasing a bit the number of corrections with respect to the theoretical ones we achieve the correct slope for the error convergence, i.e., $K\gtrsim 7$. The reason of this behaviour is still under investigation.
As we can see in figure \ref{pic:1Dscalar_convergence}, the convergence of the scheme is what we expected from theory.
\par Moreover, we can see in figure  \ref{pic:eps_vary} that, also varying the relaxation parameter $\varepsilon$, the  order of accuracy is the expected one. There are slight oscillations in particular for $\B^2$ solutions. This is a well known problem of order reduction as $\varepsilon$ is approaching the magnitude of $\Delta$, which affects lots of schemes, including some RK methods, as stated in \cite{Boscarino_Jing_Mei}. Anyway, we can say that the scheme is getting an order of accuracy bigger or equal than the expected one, except for few mid--range values of $\varepsilon$.
Moreover, we can state that the scheme is stable, for any value of $\varepsilon$ we use.
\subsubsection{Euler equation -- Isentropic flow} 
Now, we can pass to systems of equations. In particular, we will focus on Euler equation
\begin{equation}\label{eq:euler_eq}
\begin{pmatrix}
\rho\\
\rho v\\
E
\end{pmatrix}_t + \begin{pmatrix}
\rho v\\
\rho v^2 + p\\(E+p)v
\end{pmatrix}_x =0
\end{equation}
 on domain $[ -1,1 ]$, where $\rho$ is the density, $v$ the speed, $p$ the pressure and $E$ the total energy. The quantities are linked by the equation of state (EOS) \begin{equation}\label{eq:EOS_euler}
E=\frac{p}{\gamma -1} + \frac{1}{2}\rho v^2.
\end{equation}
 To test the convergence of the scheme on 1D Euler equations, we use the case of isentropic flow, when $\gamma=3$ and $p=\rho^\gamma$.
With following initial conditions
\begin{equation*}
\begin{pmatrix}
\rho_0\\
v_0\\
p_0
\end{pmatrix} = \begin{pmatrix}
1+0.5\cdot \sin(\pi x)\\
0\\ \rho_0^\gamma
\end{pmatrix}
\text{ for }x \in [-1,1 ],
\end{equation*}
final time $T=0.1$ and periodic boundary conditions.

Now, we use $\varepsilon=10^{-9}$, convection coefficient $\lambda = 3$ and CFL = $0.2$.  The $\theta$ parameter used for this convergence test, are the same of the scalar one: for $\B^1$ we used $\theta_1 = 1$, for $\B^2$ we used $\theta_1=1,\, \theta_2=0$ and for $\B^3$ we used $\theta_1=1,\,\theta_2=5$. Also here, we need a bit more of corrections for $\B^3$ to reach the $4^{\text{th}}$ order of accuracy ($K\approx 7$).
\begin{figure}[ht!]
\begin{center}
\includegraphics[scale=0.45,trim={0 15 0 25},clip]{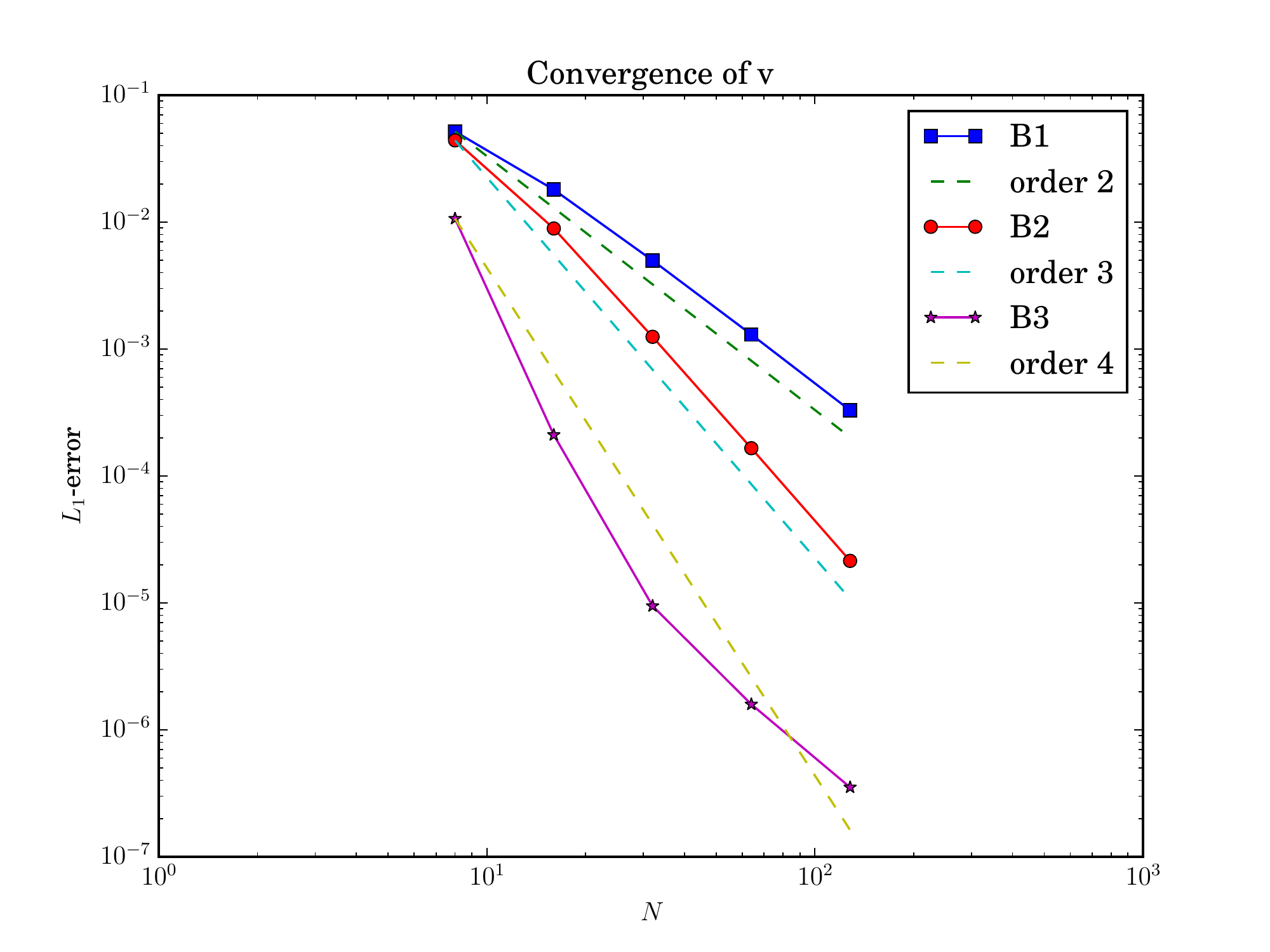}
\end{center}\caption{Convergence of Euler system}\label{pic:1Dconv_euler}
\end{figure}
As we can see in picture \ref{pic:1Dconv_euler}, the order of convergence is what we expected.
\subsubsection{Euler equation -- Sod shock test}
Now we can start testing the scheme on not smooth solutions. Let us begin with the Euler Sod test case. The Sod test case is solving equation \eqref{eq:euler_eq} on domain $[ 0,1 ]$, with EOS  $E=\frac{p}{\gamma -1} + \frac{1}{2}\rho v^2$, where $\gamma=1.4$. The initial conditions are the following
\begin{equation*}
\begin{pmatrix}
\rho_0\\
v_0\\
p_0
\end{pmatrix} = \begin{pmatrix}
1\\
0\\
1
\end{pmatrix}
\text{ for }x\leq 0.5 \qquad
\begin{pmatrix}
\rho_0\\
v_0\\
p_0
\end{pmatrix}= \begin{pmatrix}
0.125\\
0\\
0.1
\end{pmatrix}\text{ for }x>0.5,
\end{equation*}
final time is $T=0.16$ and we have outflow boundary conditions.

\begin{figure}[ht!]
\begin{center}
		\subfigure[$N=64$]{\includegraphics[width=0.49\textwidth,trim={30 20 30 20},clip=]{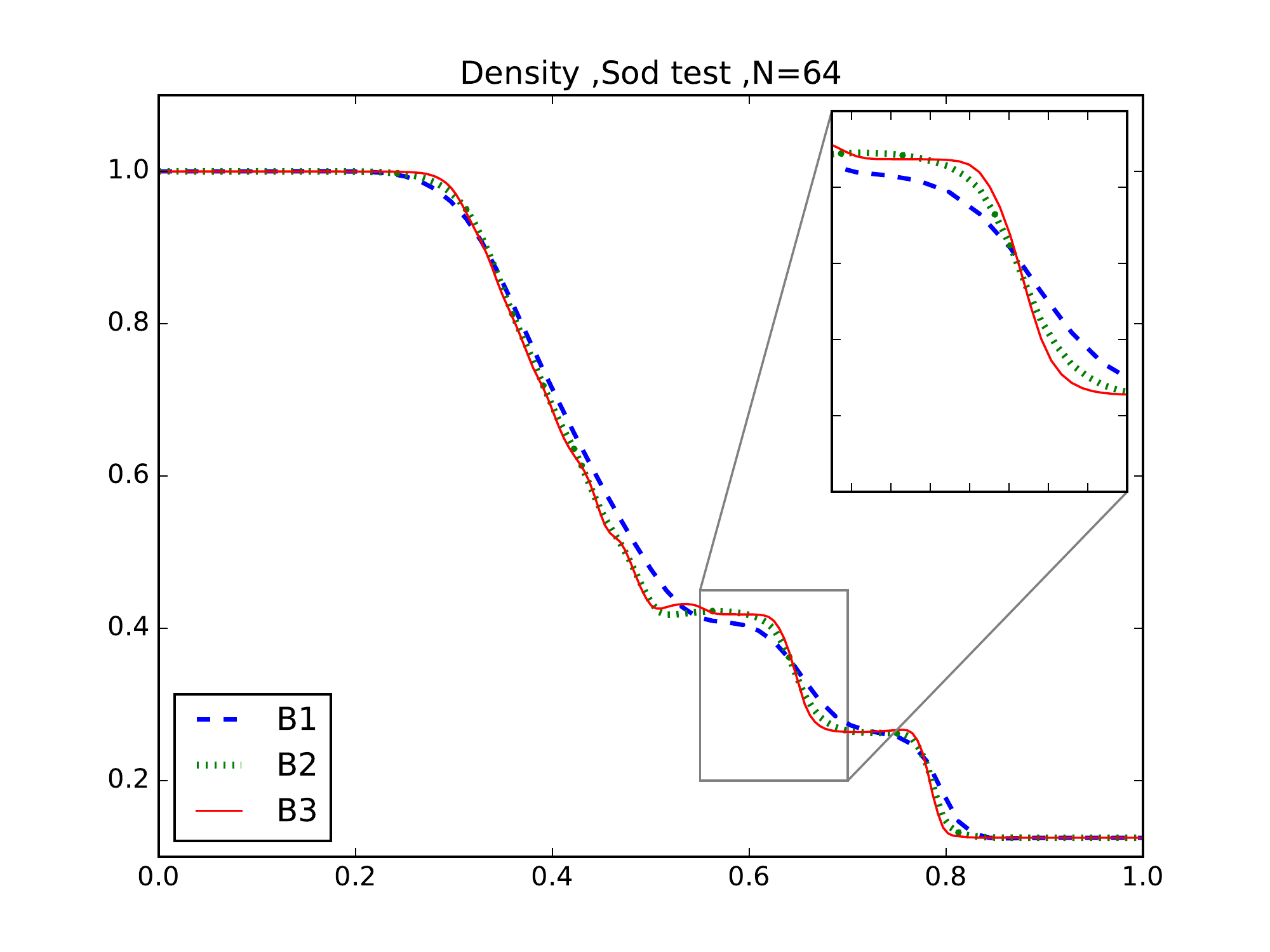}}
		\subfigure[$N=256$]{\includegraphics[width=0.49\textwidth,trim={30 20 30 20},clip=]{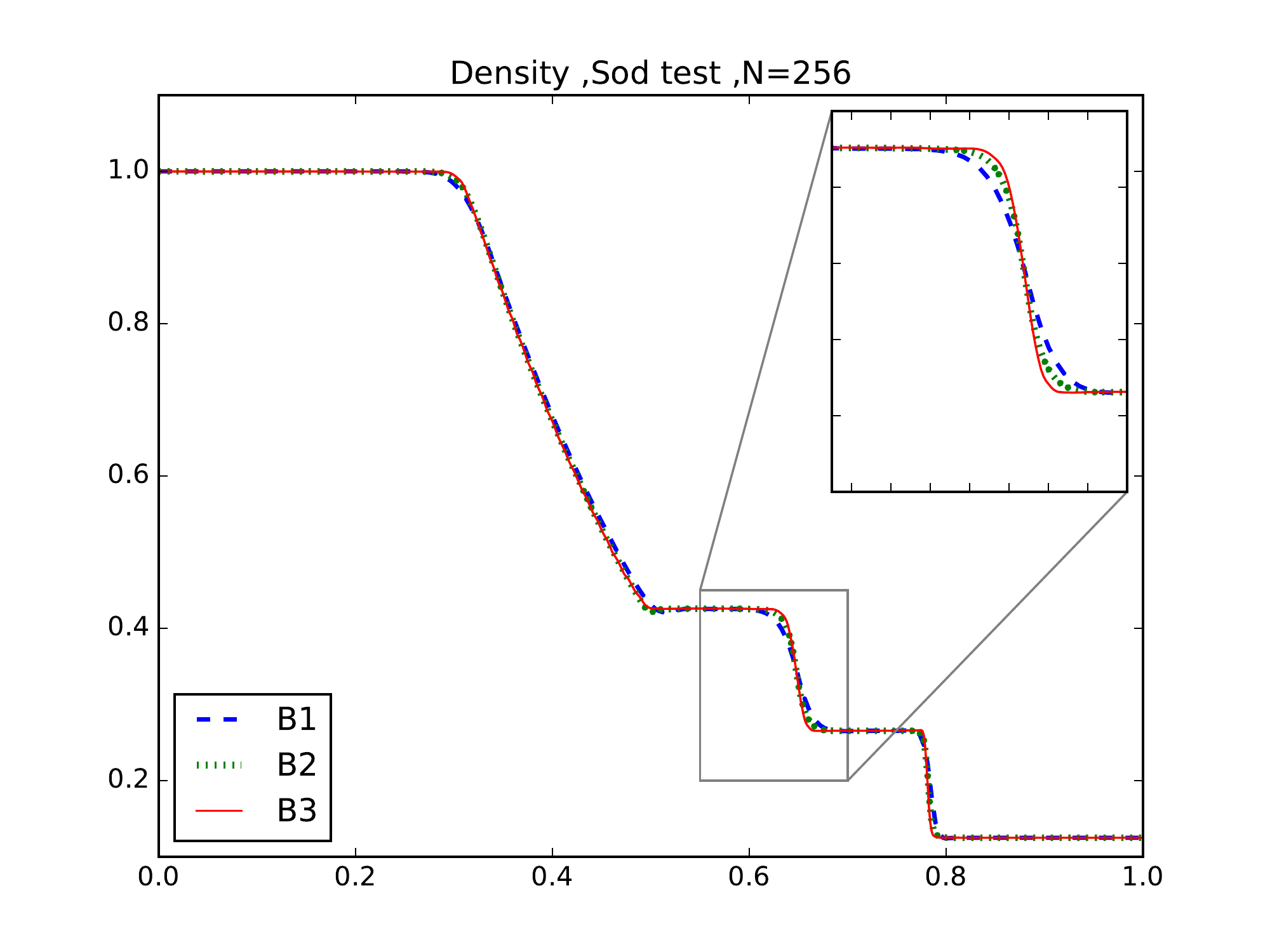}}
\end{center}\caption{Density of Sod test case 1D}\label{pic:Sod_1D}
\end{figure}

In figure \ref{pic:Sod_1D}, we can see what we obtained for $\varepsilon=10^{-9}$ in the formulation of IMEX Kinetic scheme \eqref{eq:scheme5}.
We used convection coefficient $\lambda = 2$, CFL = $0.2$. For $\B^1\; \theta_1=1$, for $\B^2\; \theta_1=1,\,\theta_2=0.5$, for $\B^3$ $ \theta_1=2.5,\,\theta_2=4$.
In picture \ref{pic:Sod_1D} we show the density plots for different mesh sizes $N=64,256$. As we can see, even with few points the $\mathbb{B}^3$ solution is outperforming the other solutions, catching in a better way the edges of the discontinuities.
\subsubsection{Euler equation -- Woodward Colella}
\begin{figure}[ht!]
\begin{center}
\subfigure[$N=256$]{\includegraphics[width=0.49\textwidth,trim={0 20 0 20},clip]{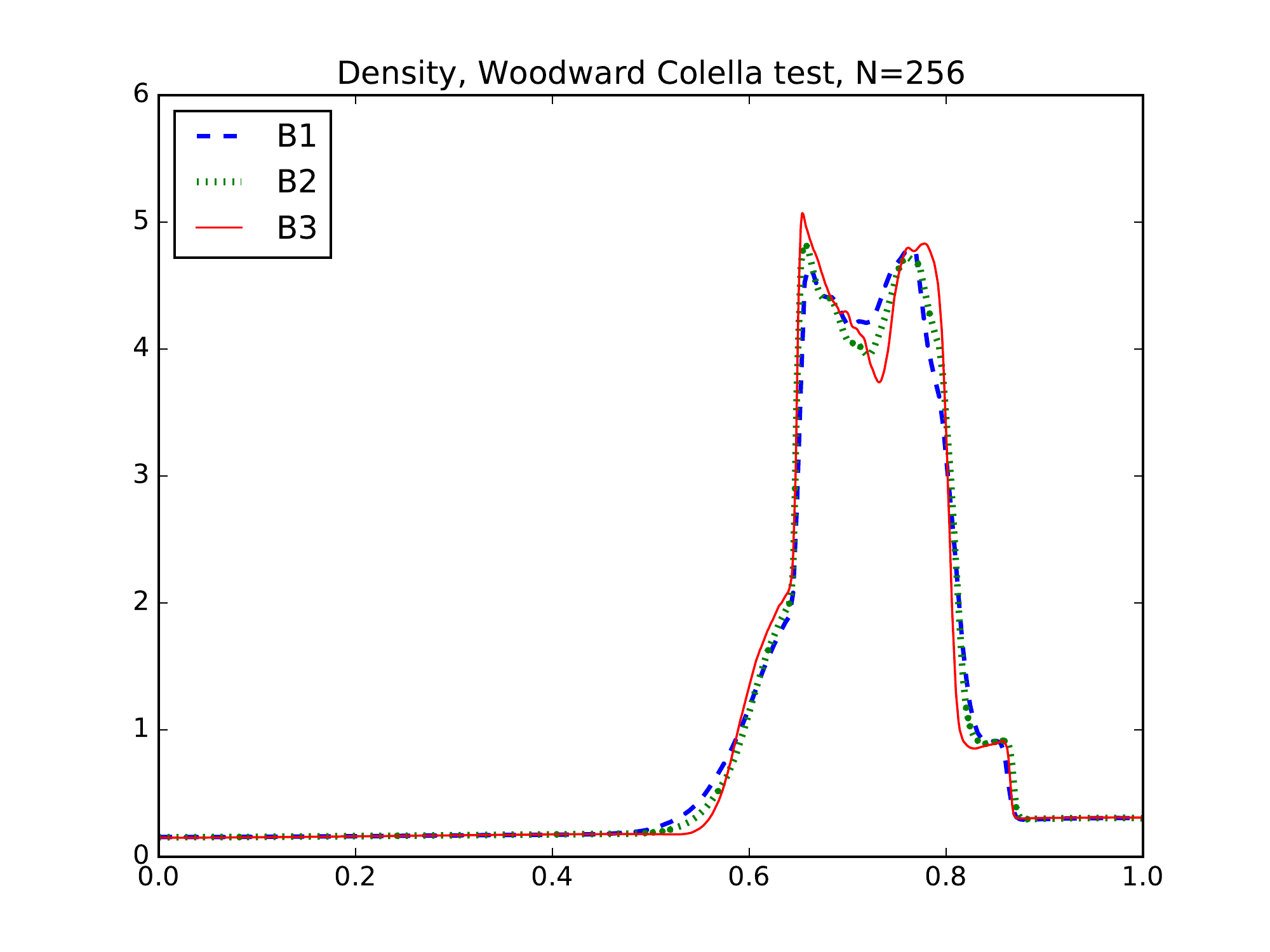}}
\subfigure[$N=512$]{\includegraphics[width=0.49\textwidth,trim={0 20 0 20},clip]{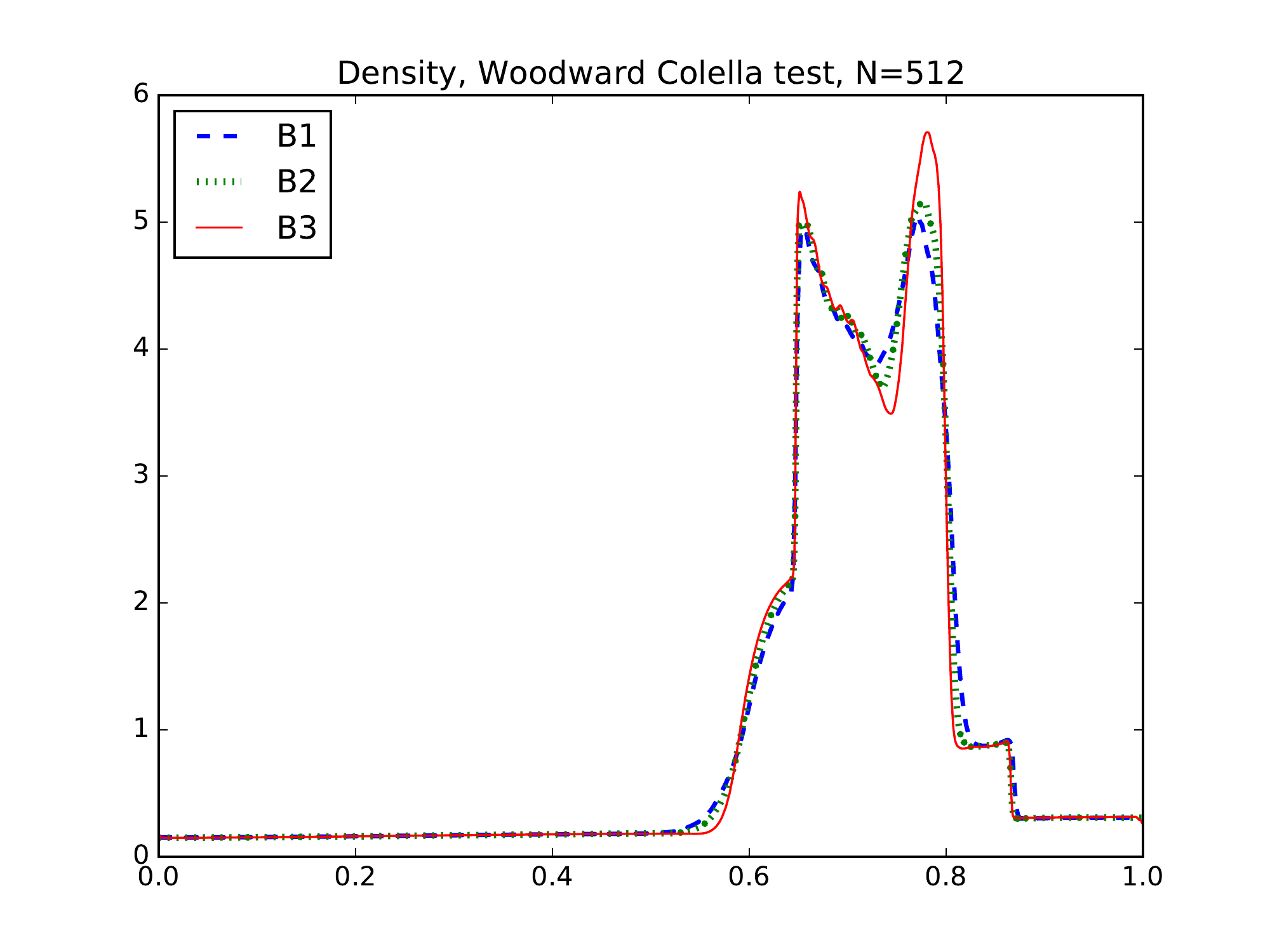}}
\end{center}\caption{Density of Woodward Colella test}\label{pic:wood}
\end{figure}
We can see even better the advantages of using a high order scheme in the following examples. First, we present the one proposed by Woodward and Colella \cite{woodward_colella}. It solves again Euler equation  \eqref{eq:euler_eq} on domain $[0,1]$ with EOS \eqref{eq:EOS_euler} with $\gamma=1.4$. The final time is 0.038, the initial conditions are 
\begin{equation*}
\rho_0=1, \qquad v_0=0, \qquad  p_0 =\begin{cases} 
10^{3} \text{ for } x\in [0,0.1],  \\
10^{-2} \text{ for } x\in [0.1,0.9],\\
10^{2} \text{ for } x\in[0.9,1]
\end{cases} 
\end{equation*}
and we use outflow boundary conditions. For $\B^1$,  $\theta_1=0.5$.  For $\B^2$, $ \theta_1=0.8$, $\theta_2=1$. For $\B^3$,  $\theta_1=5$, $\theta_2=1$. In figure \ref{pic:wood} there is the result for $\varepsilon=10^{-9}$, convection coefficient = $20$, CFL = $0.1$, $N=256,512$. 

We can notice that in this case, only $\mathbb{B}^3$ is able to catch the shape of the second peak (with 512 elements).

\subsubsection{Euler equation -- Shu Osher test}
Last test we performed in 1D was proposed by Shu and Osher \cite{shu_osher}. Again we have Euler equation \eqref{eq:euler_eq}  on domain $[ -5,5 ]$ with EOS \eqref{eq:EOS_euler} with $\gamma=1.4$. Here initial conditions are
\begin{equation*}
\begin{pmatrix}
\rho_0\\
v_0\\
p_0
\end{pmatrix} = \begin{pmatrix}
3.857143\\
2.629369\\
10.333333
\end{pmatrix}
\text{ if } x\in [-5, -4], \quad
\begin{pmatrix}
\rho_0\\
v_0\\
p_0
\end{pmatrix}= \begin{pmatrix}
1 + 0.2\sin(5x)\\
0\\
1
\end{pmatrix}\text{ if } x\in[-4, 5].
\end{equation*}
\begin{figure}[ht!]
\begin{center}
	\subfigure[$N=64$]{\includegraphics[width=0.49\textwidth,trim={0 20 0 20},clip]{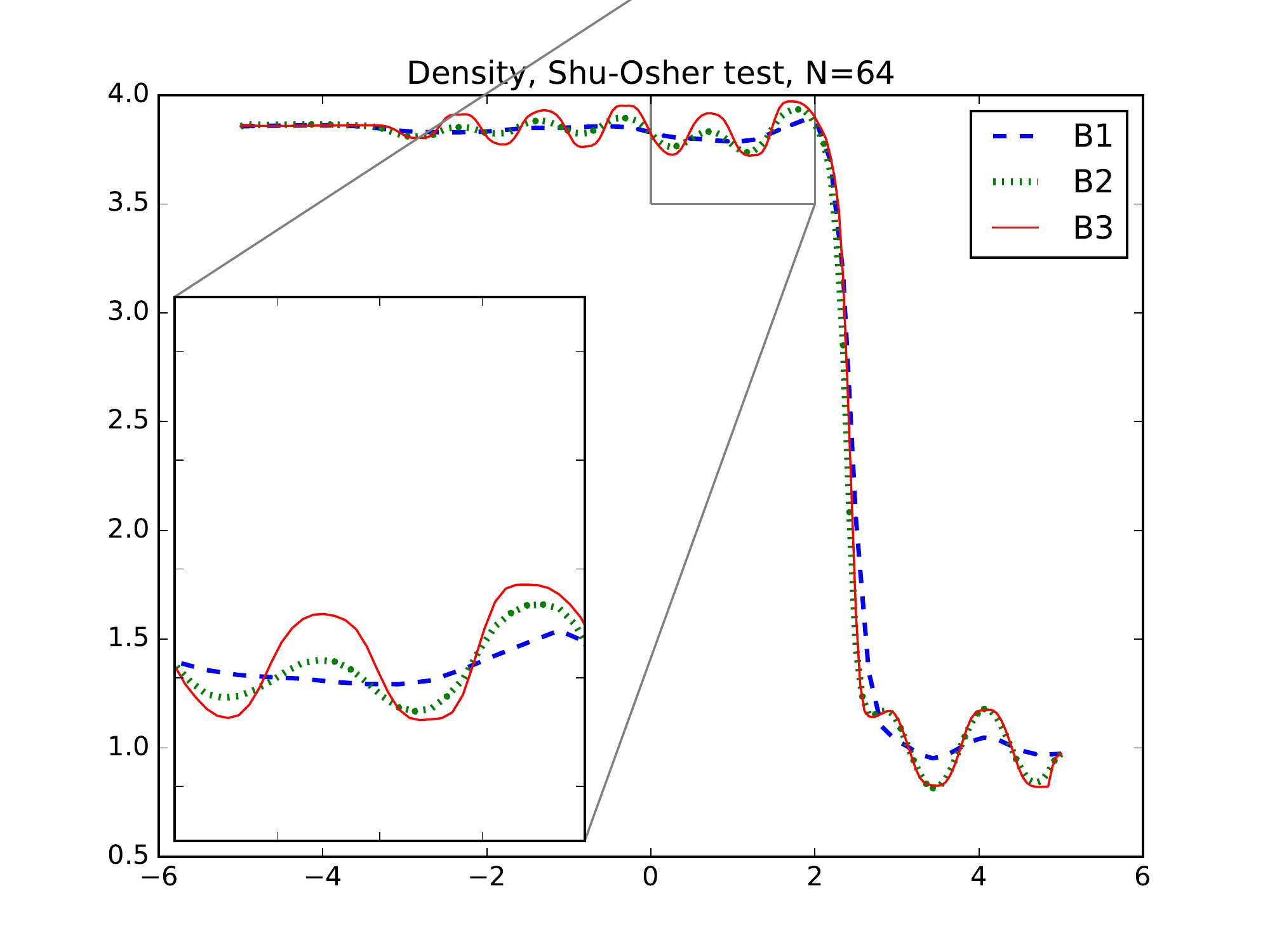}}
	\subfigure[$N=128$]{\includegraphics[width=0.49\textwidth,trim={0 20 0 20},clip]{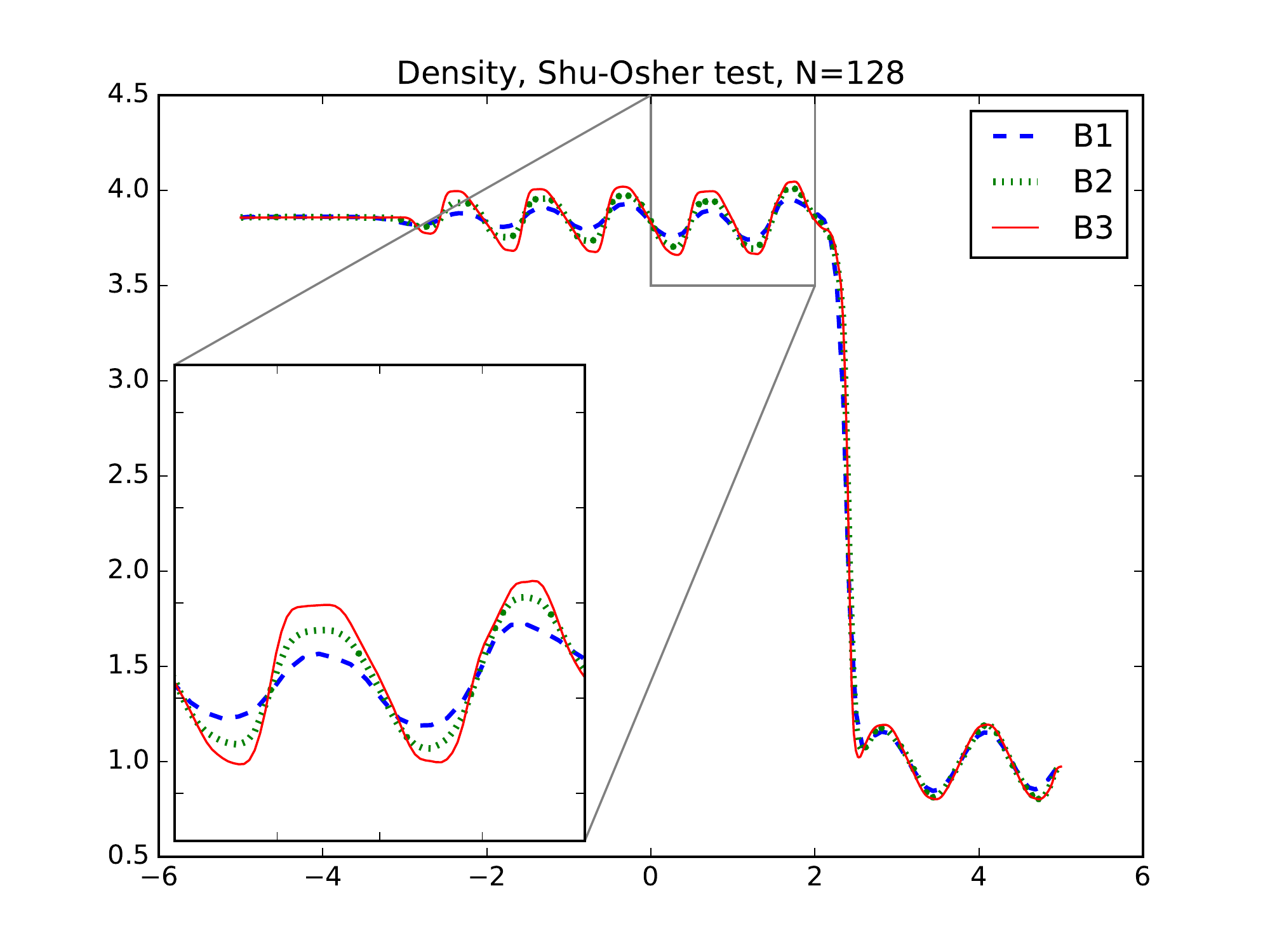}}\\
	\subfigure[$N=256$]{\includegraphics[width=0.49\textwidth,trim={0 20 0 20},clip]{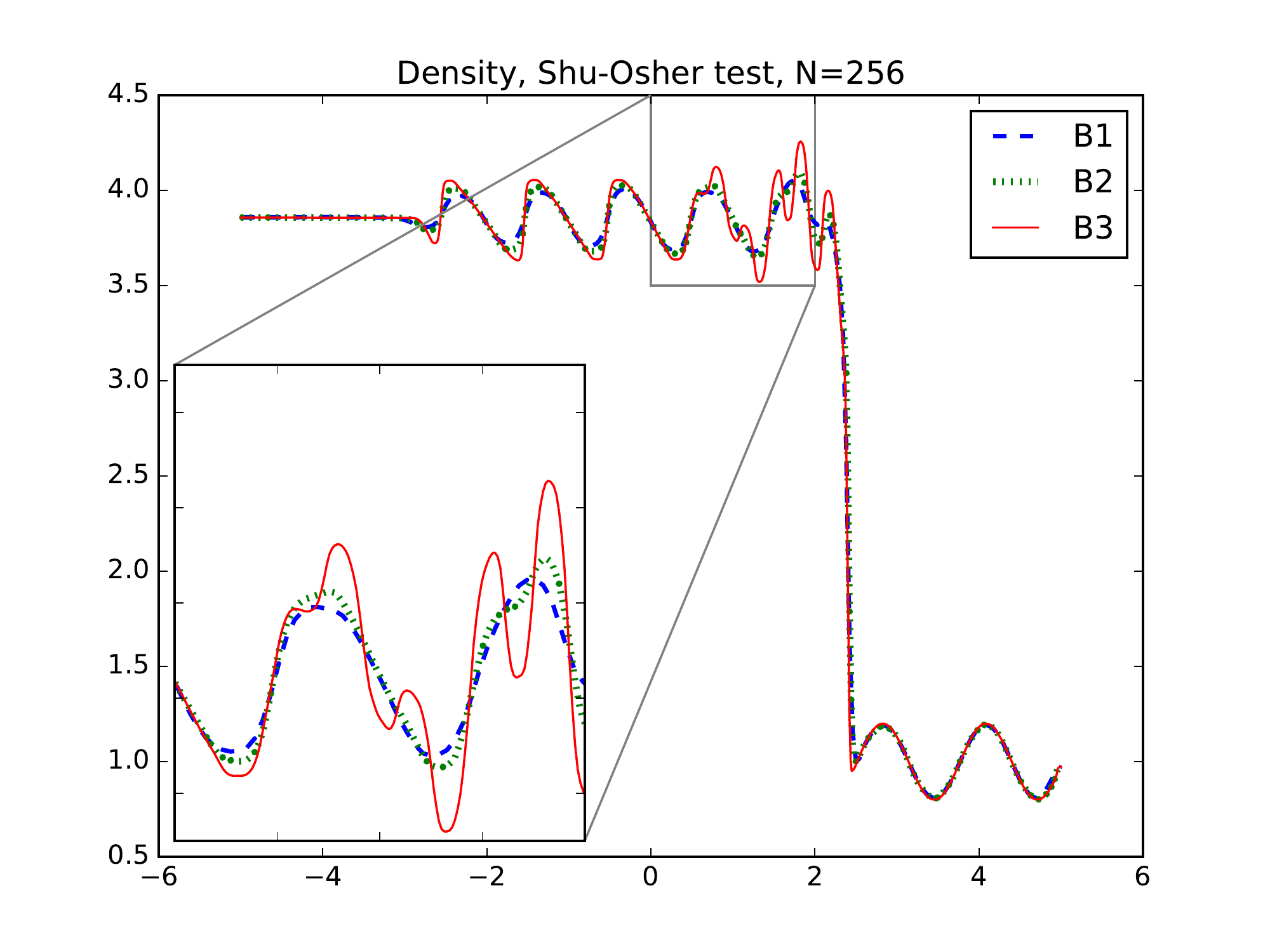}}
	\subfigure[$N=512$]{\includegraphics[width=0.49\textwidth,trim={0 20 0 20},clip]{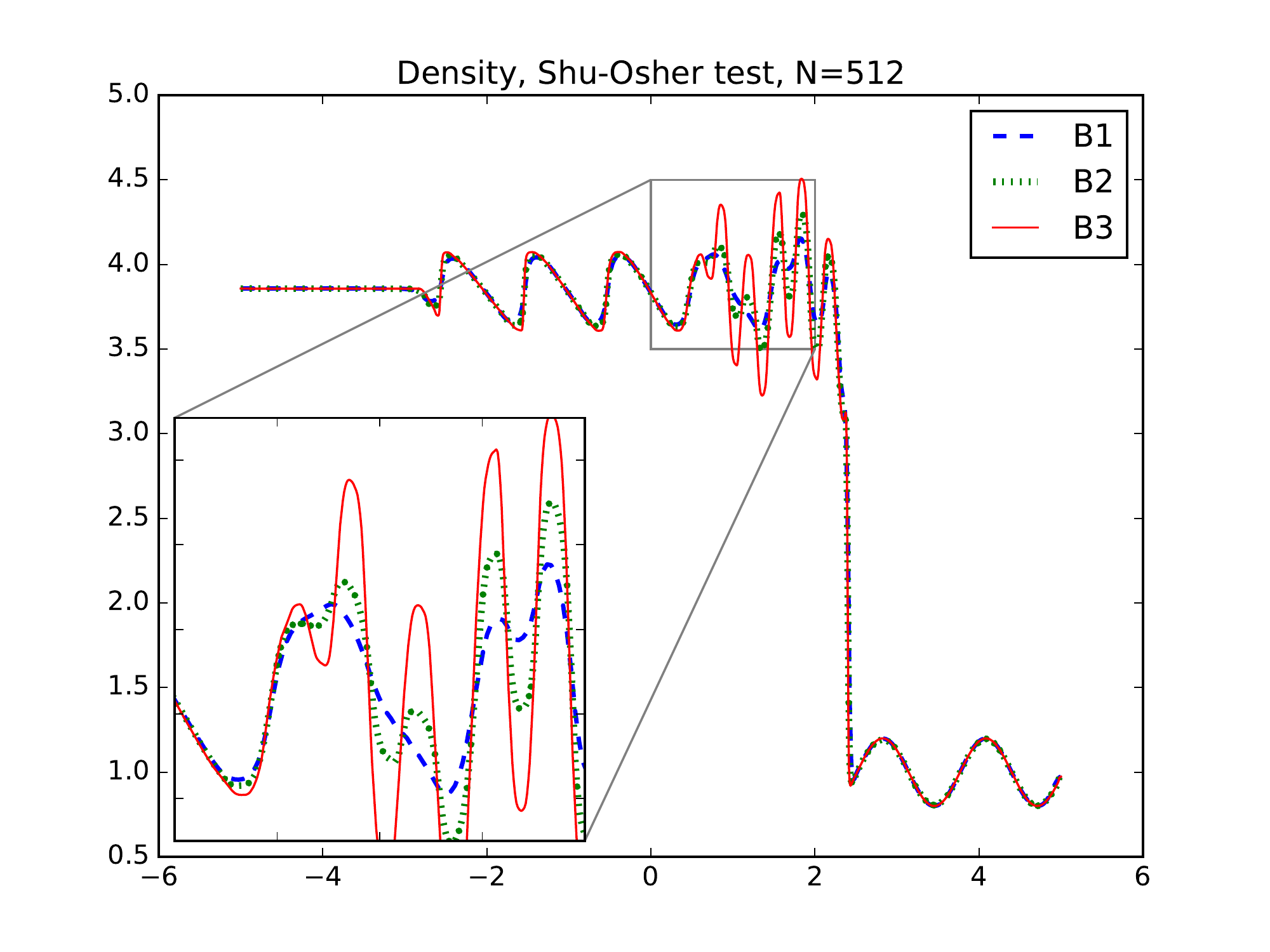}}
\end{center}\caption{Density of Shu--Osher's test }\label{pic:shu}
\end{figure}

Final time is $T=1.8$, we use outflow boundary conditions, $\varepsilon=10^{-9}$, convection coefficient $\lambda=3$, CFL = $0.1$.   For $\B^1\, \theta_1=0.5$, for $\B^2\, \theta_1=0.8,\,\theta_2=1$, for $\B^3\, \theta_1=3,\,\theta_2=1$. In figure  \ref{pic:shu}, we can see results for several $N$s. Even here, we can see that the second and third order polynomials perform better with respect to the first order one. In particular, we can see how the oscillations are already captured with few points and how the precision increases quickly if the order is greater.

In all these cases, we have seen that our method performs nicely and capture the correct behaviours of the equations solutions. Moreover, we see that it can be convenient to switch to higher order to better get the solution of our test cases with less mesh elements.  

\subsection{2D numerical tests}
Let us present some numerical test defined on a 2D domain. We use again the DRM model poposed by \cite{natalini} and the scheme we presented. We see only examples of Euler equation in 2D:
\begin{equation}\label{eq:euler_2D}
\begin{split}
\partial_t U(\mathbf{x},t) + \partial_{x} A_1(U(\mathbf{x},t)) +\partial_y A_2(U(\mathbf{x},t))=0,\qquad \mathbf{x} =(x,y) \in \Omega\subset \R^2, \\
U=\begin{pmatrix}
\rho \\ \rho u \\ \rho v \\ E
\end{pmatrix}, \qquad A_1 (U)=\begin{pmatrix}
\rho u \\ \rho u^2+p \\ \rho u v \\  u(E+p)	
\end{pmatrix}
, \qquad A_2(U)=\begin{pmatrix}
\rho v \\ \rho uv \\ \rho v^2 +p \\  v(E+p)	
\end{pmatrix}
\end{split}
\end{equation}
where $\rho$ is the density, $u$ is the speed in $x$ direction, $v$ is the speed in $y$ direction, $E$ the total energy and $p$ the pressure. They are linked by the following EOS:
\begin{equation}\label{eq:EOS_2D}
p=(\gamma-1)\Big( E-\frac{1}{2}\rho (u^2+v^2)\Big).
\end{equation}
\subsubsection{Euler equation -- Smooth vortex test case}
\begin{figure}[ht!]
\begin{center}
\includegraphics[scale=0.5,trim={20 10 20 20},clip]{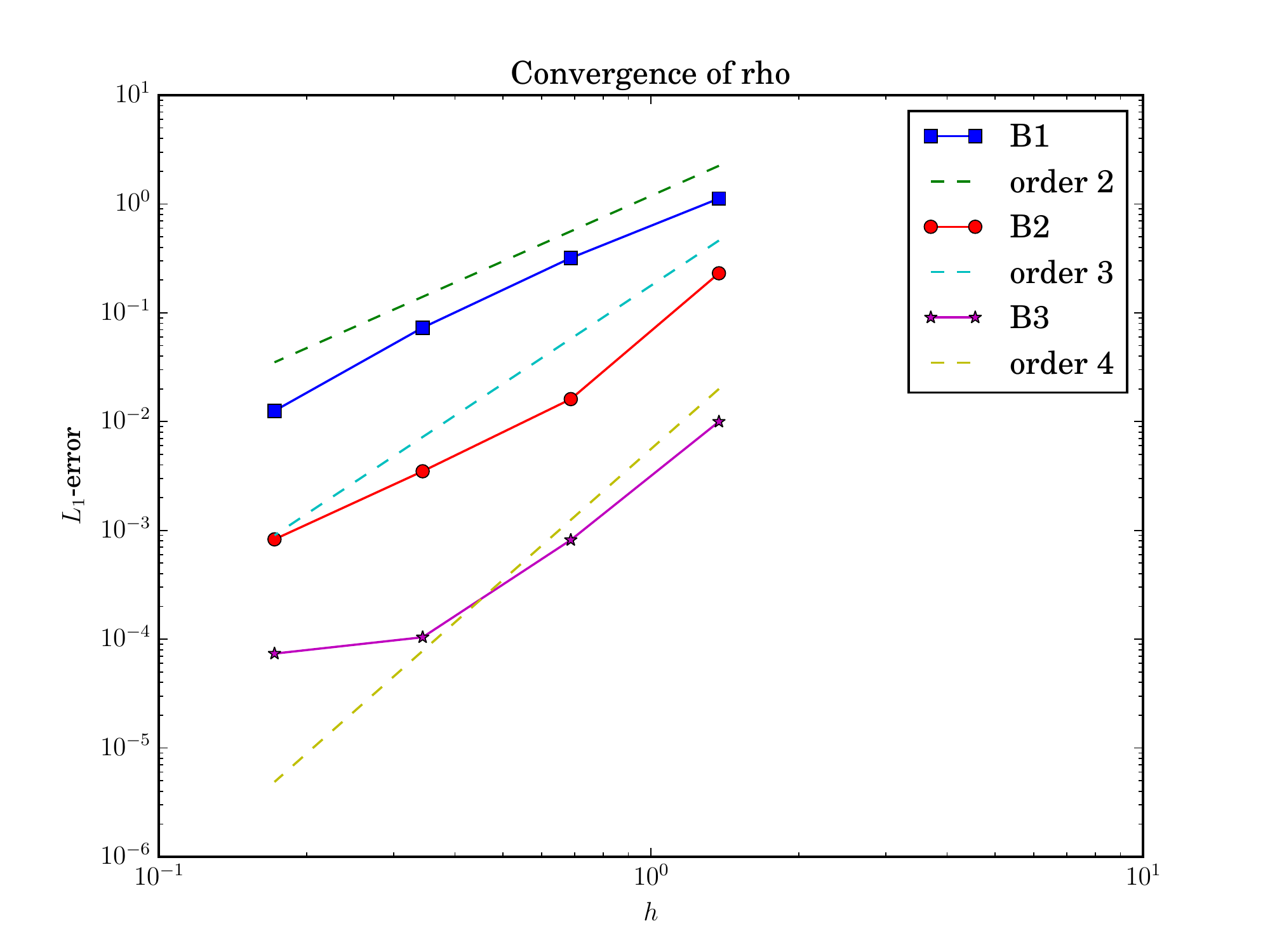}
\end{center}\caption{2D convergence} \label{pic:2D_conv}
\end{figure}
To start, we want to study the convergence of the method also in 2D. To do so, we test our scheme with a steady vortex test case, so that we can compare the final solution with the initial one. The domain is a circle of radius 10 and center $(0,0)$. The initial conditions are
$$
\begin{pmatrix}
\rho_0\\ u_0\\ v_0 \\ p_0
\end{pmatrix} =
\begin{pmatrix}
\left ( 1-\frac{\gamma-1}{\gamma}\frac{1}{2} \left(\frac{5}{2\pi}\right)^2  e^{\frac{1-r^2}{2}} \right)^\frac{1}{\gamma-1}\\
\frac{5}{2\pi} (-y)  e^{\frac{1-r^2}{2}}\\
	\frac{5}{2\pi} (x)  e^{\frac{1-r^2}{2}}\\
		\rho_0^\gamma
\end{pmatrix}.
$$
Here $r^2=x^2+y^2$ and the boundary conditions are outflow.
In our simulations $\gamma = 1.4$ for the EOS \eqref{eq:EOS_2D}. Again, we take $\varepsilon=10^{-9}$, convection coefficient $\lambda =1.4$ and CFL = $0.1$. We stop the simulation at time $T=1.$ We use different refinements of the domain mesh. These are uniform triangular meshes and on the x--axis of figure \ref{pic:2D_conv} one can see the maximum diameter of a cell of the mesh. We can see in figure \ref{pic:2D_conv} that the convergence is reflecting the theoretical results, even if for $\B^3$ we need more corrections ($K\approx 7$) to get the order to get closer to the convergence expected. For $\B^1\, \theta_1=0.1$, for $\B^2\, \theta_1=0.01,\,\theta_2=0$, for $\B^3\, \theta_1=0.001,\,\theta_2=0$.

\subsubsection{Euler equation -- Sod 2D test case}
We tested our method on the analogous of Sod in 2D. This test is again solving Euler equation \eqref{eq:euler_2D} where $\gamma=1.4$ in EOS \eqref{eq:EOS_2D}. The domain $\Omega$ is a circle of radius 1 and center in $(0,0)$. The initial conditions are:
$$
\begin{pmatrix}
\rho_0\\ u_0\\ v_0 \\ p_0
\end{pmatrix} =
\begin{pmatrix}
1\\
 0 \\ 0 \\ 1
\end{pmatrix} \text{ if } x^2+y^2<\frac{1}{4}, \qquad \begin{pmatrix}
\rho_0\\ u_0\\ v_0 \\ p_0
\end{pmatrix} =
\begin{pmatrix}
0.125\\
0 \\ 0 \\ 0.1
\end{pmatrix} \text{ if } x^2+y^2\geq\frac{1}{4}.
$$
\begin{figure}[ht!]
	\begin{center}
		\subfigure[$\mathbb{B}^1, N=3576$]{\includegraphics[width=0.49\textwidth,trim={0 100 0 90},clip]{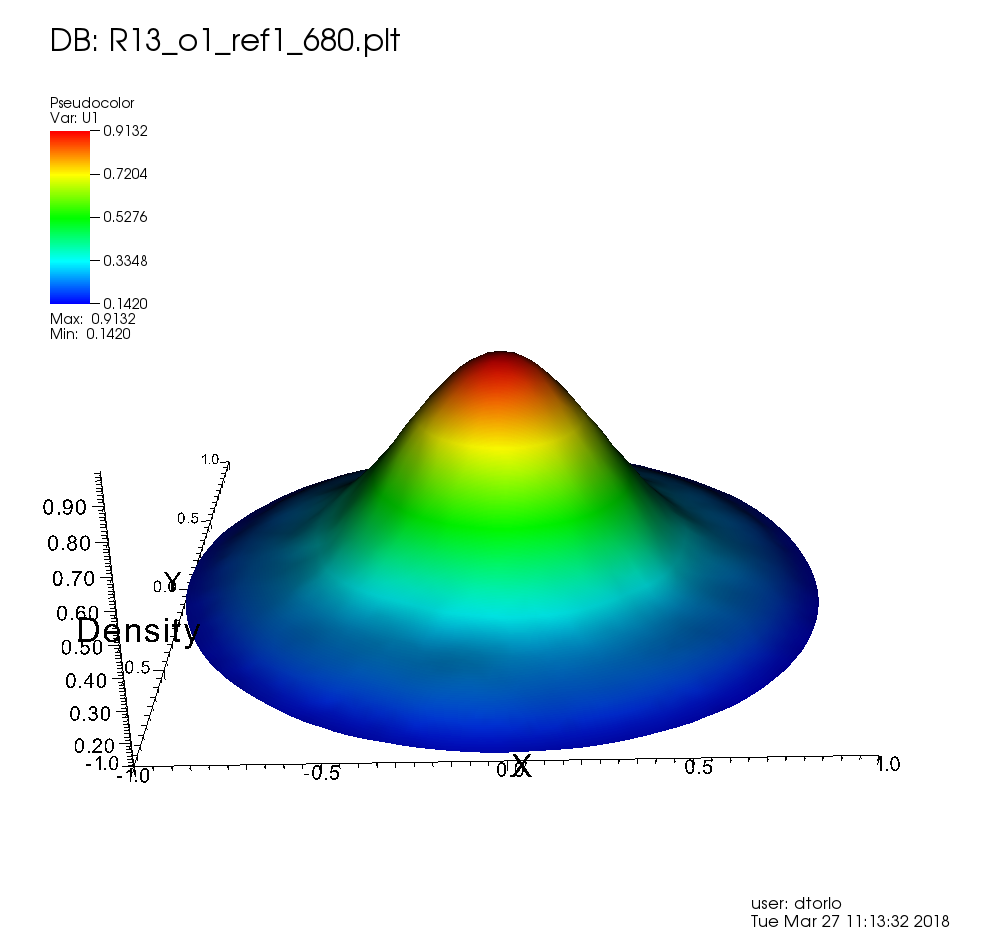}}
		\subfigure[$\mathbb{B}^1, N=13548$]{\includegraphics[width=0.49\textwidth,trim={0 100 0 90},clip]{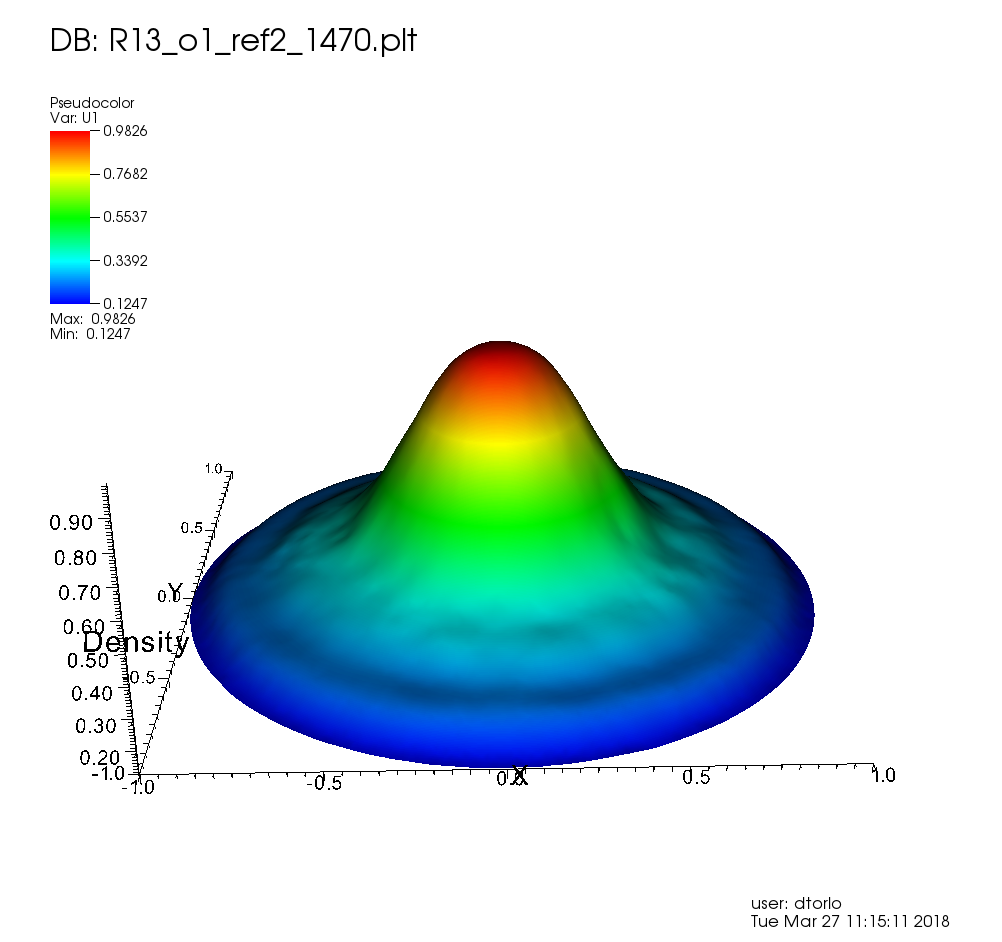}}\\
		\subfigure[$\mathbb{B}^2, N=3576$]{\includegraphics[width=0.49\textwidth,trim={0 100 0 90},clip]{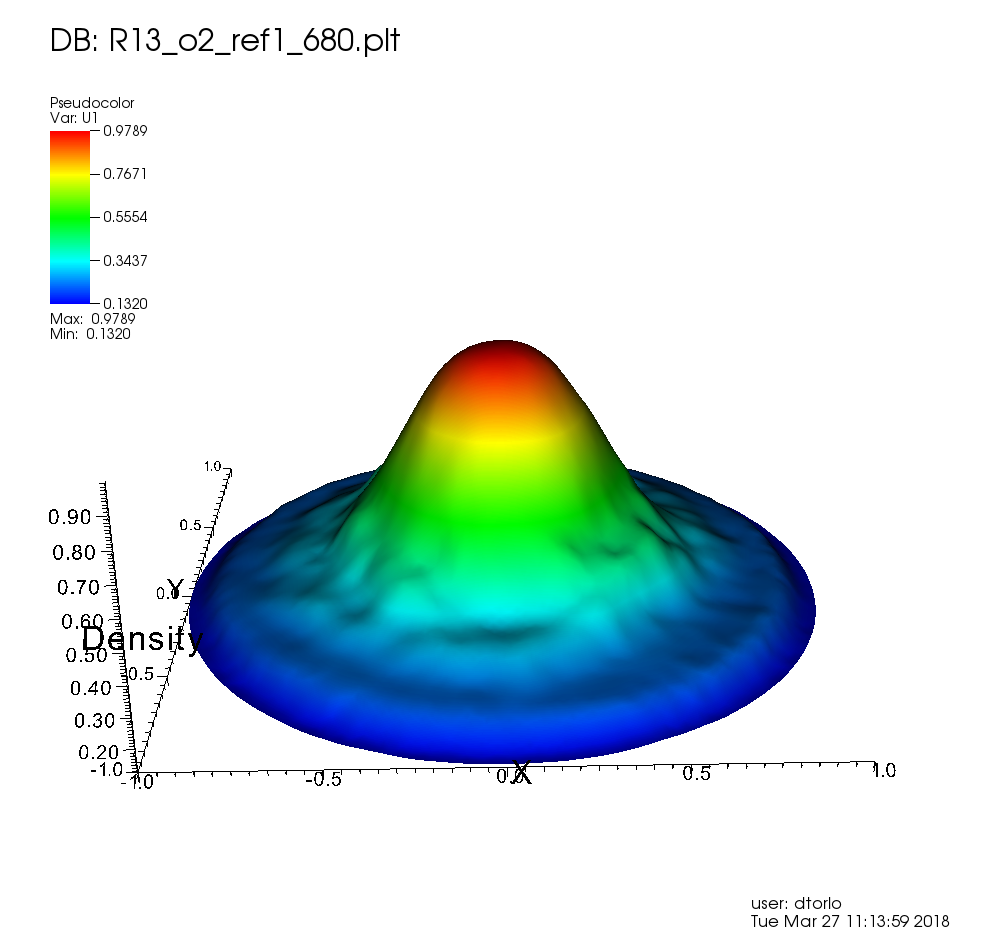}}
		\subfigure[$\mathbb{B}^2, N=13548$]{\includegraphics[width=0.49\textwidth,trim={0 100 0 90},clip]{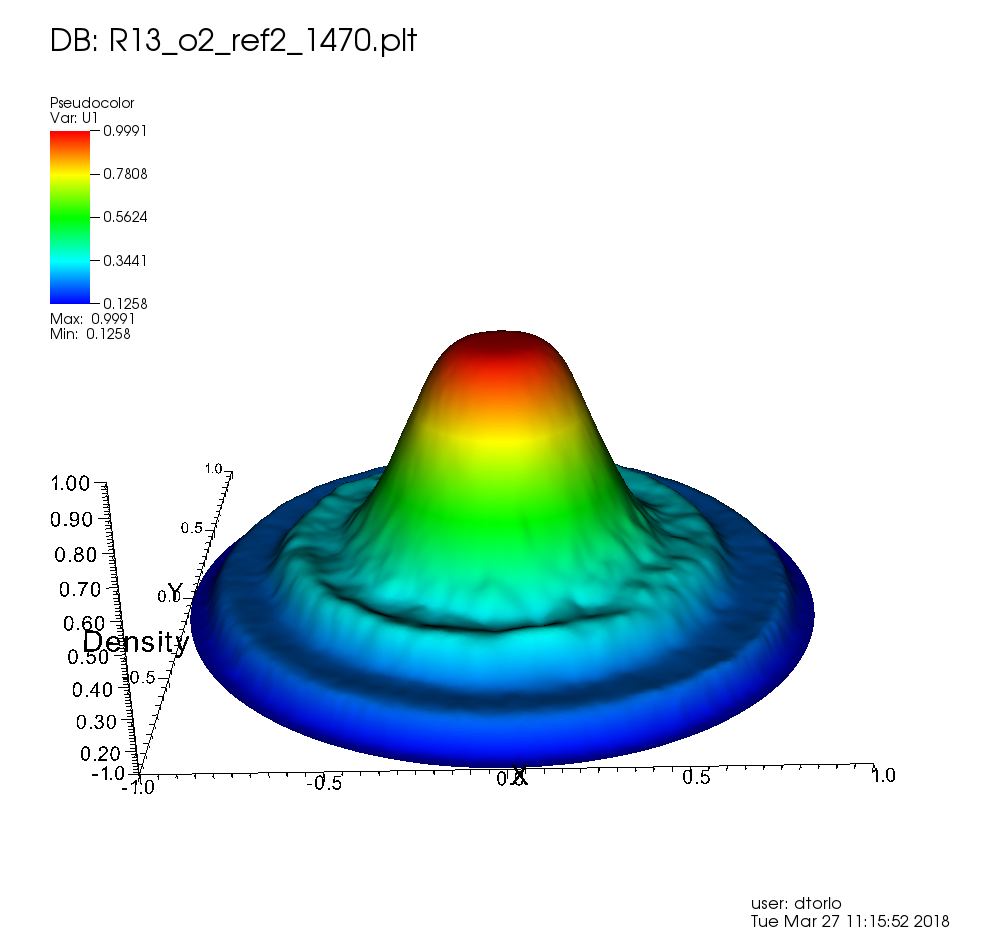}}\\
		\subfigure[$\mathbb{B}^3, N=3576$]{\includegraphics[width=0.49\textwidth,trim={0 100 0 90},clip]{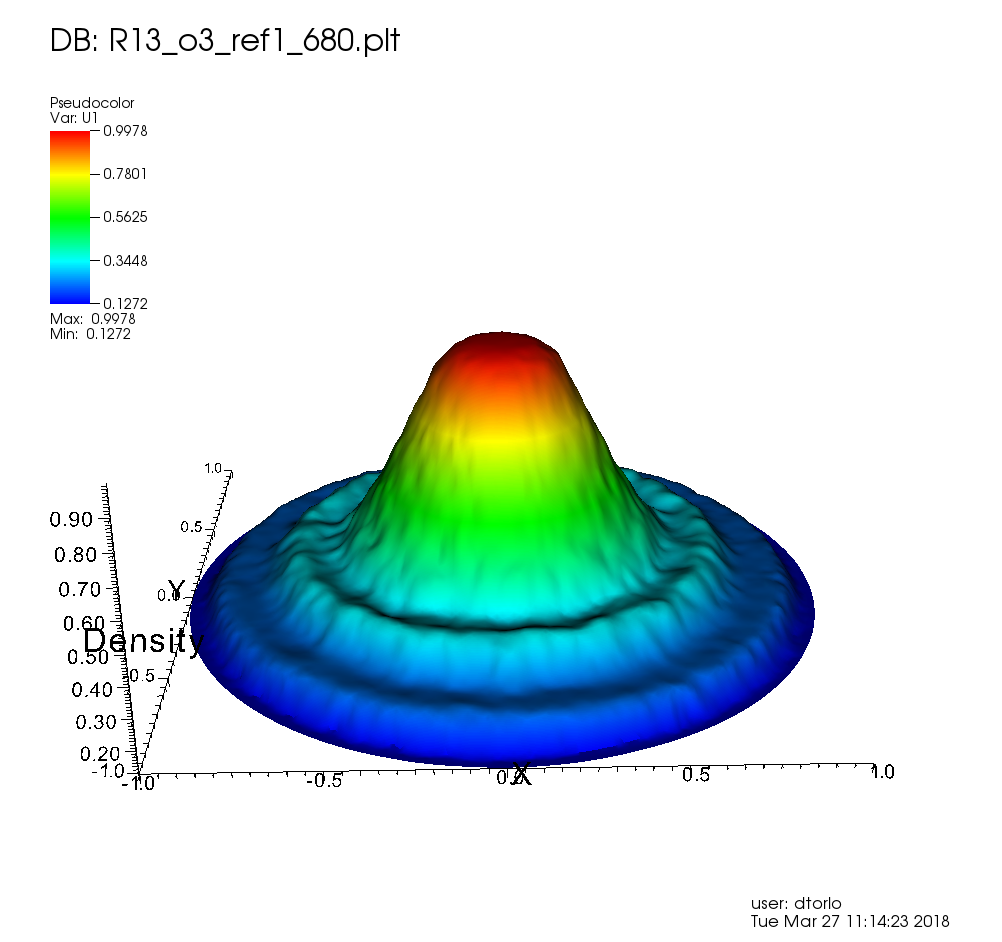}}
		\subfigure[$\mathbb{B}^3, N=13548$]{\includegraphics[width=0.49\textwidth,trim={0 100 0 90},clip]{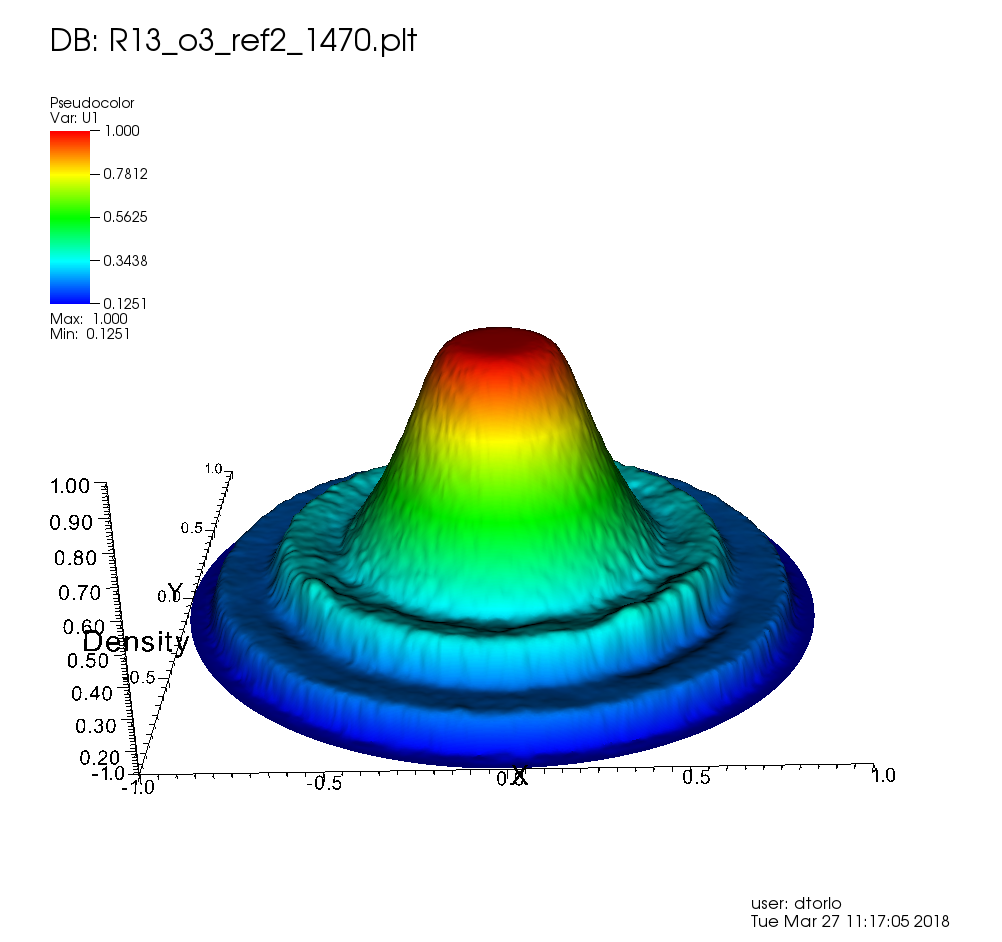}}
	\end{center}\caption{Density of Sod test}\label{pic:sod_2D_surf}
\end{figure}

\begin{figure}[ht!]
	\begin{center}
		\subfigure[Slice of $\mathbb{B}^k, N=3576$]{\includegraphics[width=0.49\textwidth,trim={140 100 0 450},clip]{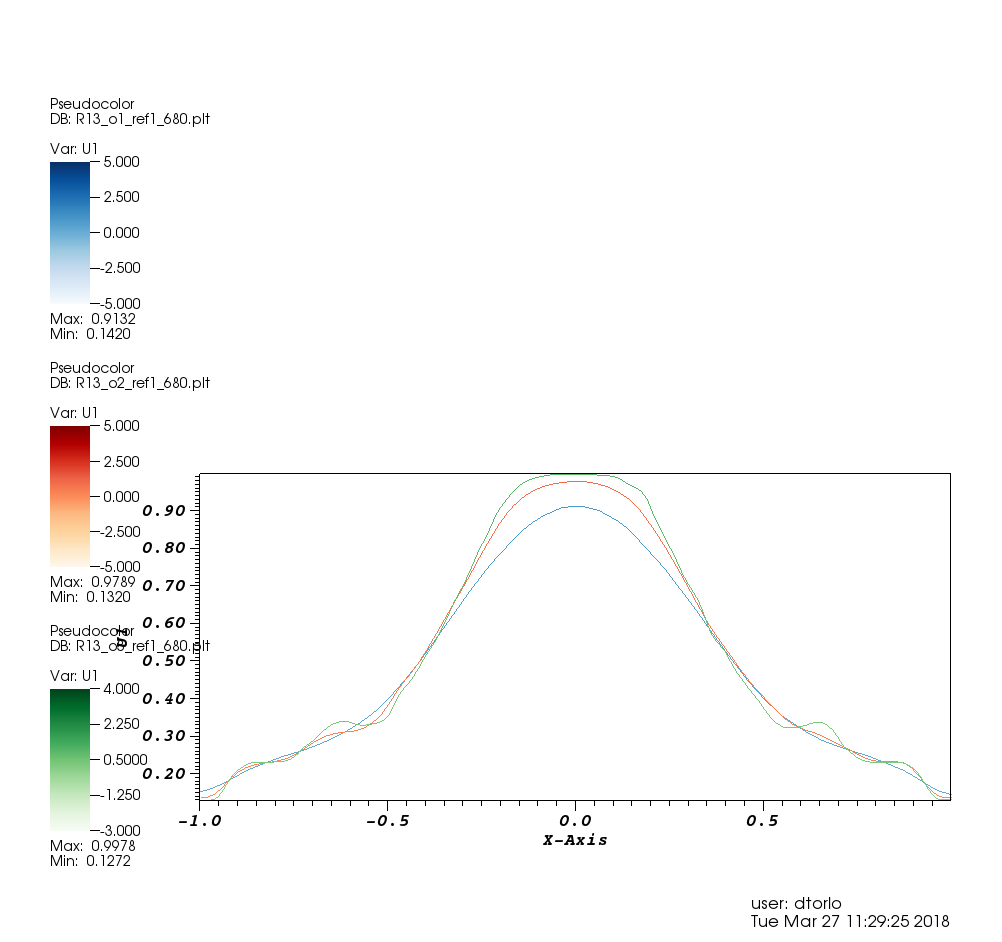}}
		\subfigure[Slice of $\mathbb{B}^k, N=13548$]{\includegraphics[width=0.49\textwidth,trim={140 100 0 450},clip]{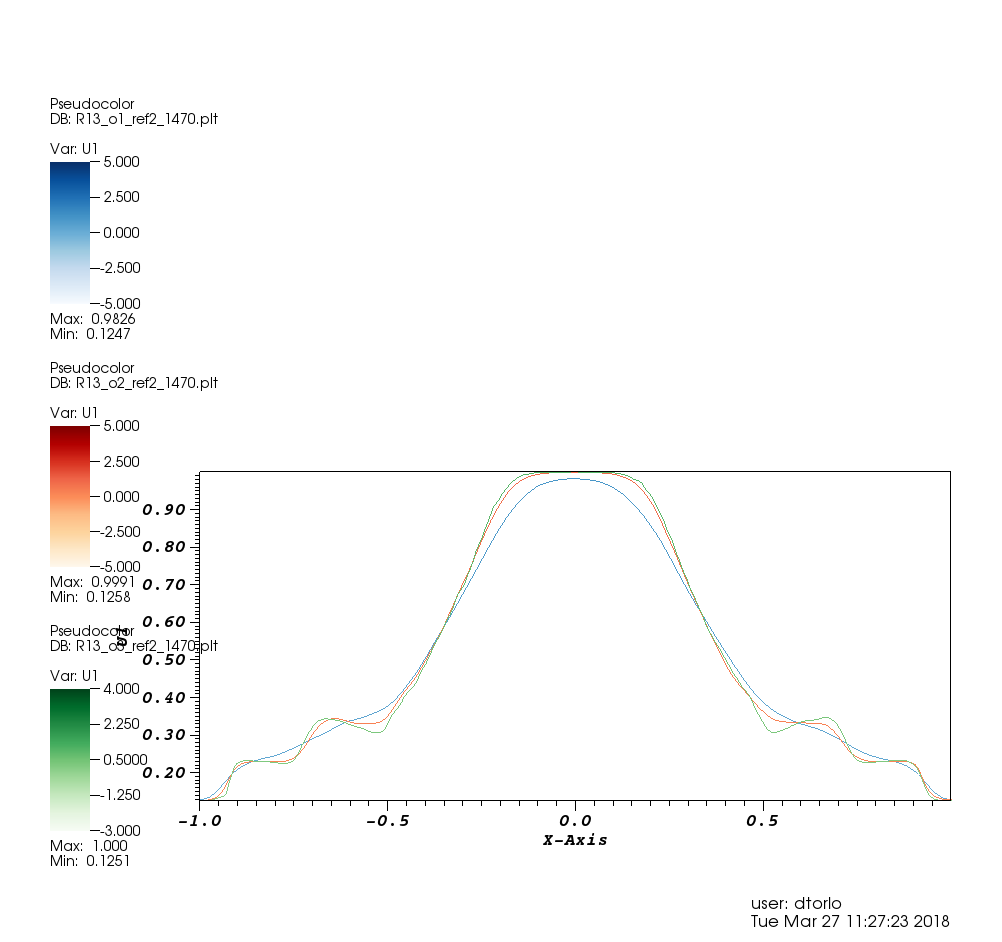}}\\
		\subfigure[Scatter of $\mathbb{B}^k, N=3576$]{\includegraphics[width=0.49\textwidth,trim={0 100 0 90},clip]{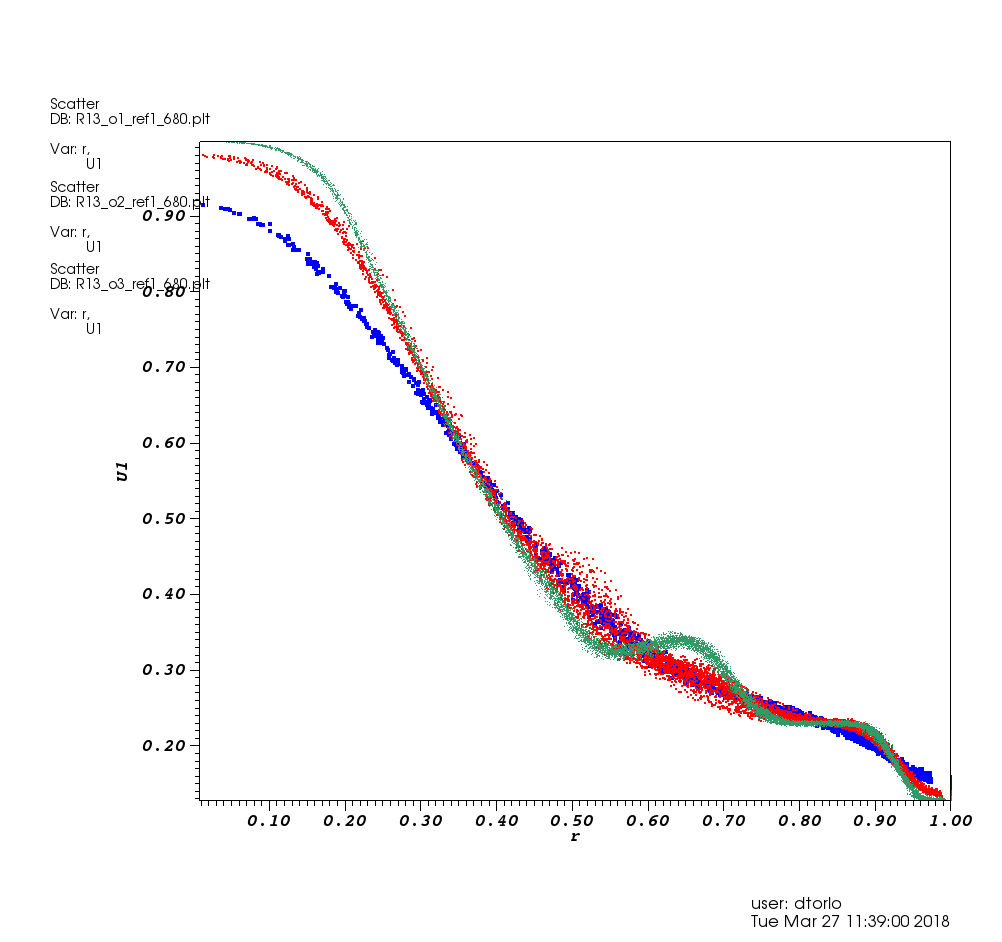}}
		\subfigure[Scatter of $\mathbb{B}^k, N=13548$]{\includegraphics[width=0.49\textwidth,trim={0 100 0 90},clip]{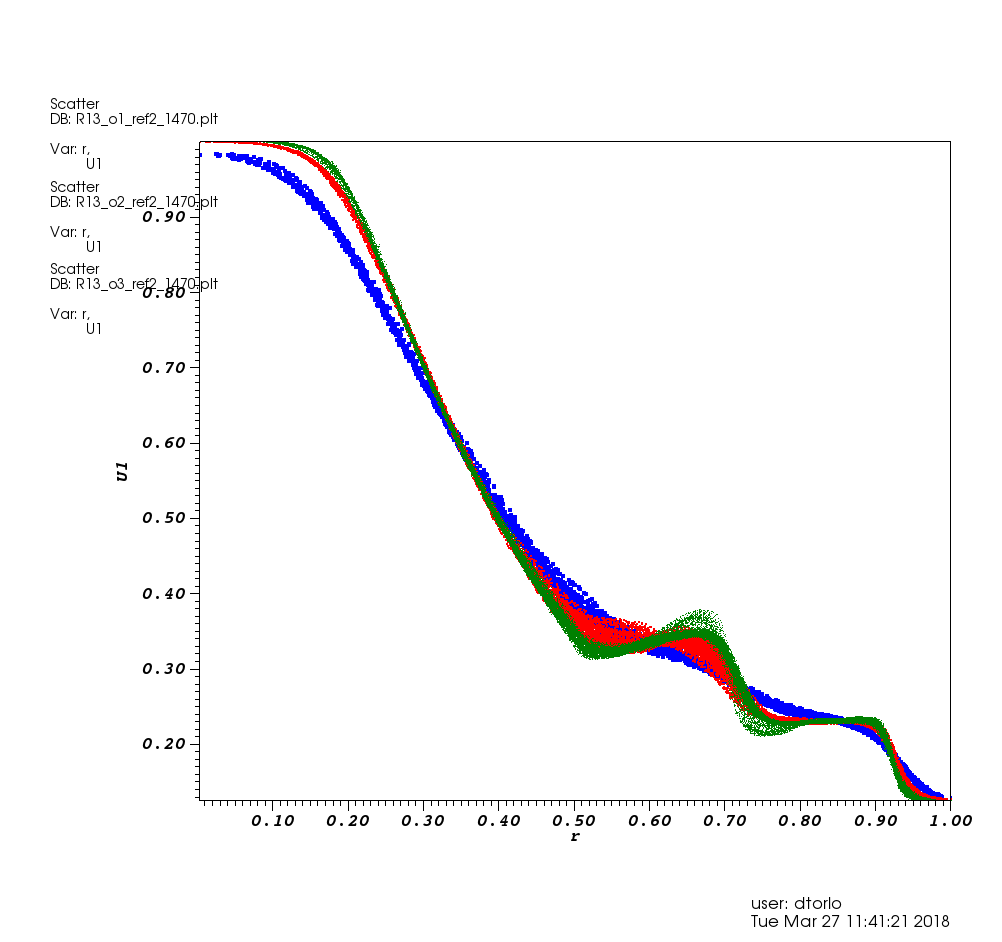}}
	\end{center}\caption{Density of Sod test ($\mathbb{B}^1$ blue, $\mathbb{B}^2$ red and $\mathbb{B}^3$ green)}\label{pic:sod_2D_slice}
\end{figure}

The parameters used for this test are $\varepsilon=10^{-9}$, convection coefficient $\lambda=1.4$, CFL = $0.1$, final time $T=0.25$ and outflow boundary conditions.  For $\B^1\, \theta_1=0.1$, for $\B^2\, \theta_1=0.1,\,\theta_2=0.0001$, for $\B^3\, \theta_1=0.01,\,\theta_2=0.0001$.

We use uniform triangular meshes with different sizes and what is shown in pictures \ref{pic:sod_2D_surf} is obtained with $N=3576$ and $N=13548$ triangles on the domain.

If we watch pictures \ref{pic:sod_2D_surf} and \ref{pic:sod_2D_slice}, we can see that also in this case the higher the order of polynomial we use, the sharper becomes the solution. In particular, we can say that the solution with $\mathbb{B}^2$ basis functions for the mesh with $N=13548$ elements is comparable with the solution for $\mathbb{B}^3$ with only $N=3576$ elements. Moreover, we can see that with $\mathbb{B}^1$ the diffusion is too high and it is smoothening all the discontinuities.

\subsubsection{Euler equation -- DMR 2D test case}
For the last test case, we test our scheme on the DMR (double Mach reflection) problem presented in \cite{DMR}. The equation we are solving is again the Euler equation \eqref{eq:euler_2D} with $\gamma =1.4$ in EOS \eqref{eq:EOS_2D}. The domain is a rectangular shape, cut on the bottom right part by an oblique edge. The boundaries of the rectangle are $x=0,\, x=2.2,\,y=-0.2, \, y=3$. The oblique edge is a line passing through points $(0,0)$ and $(3,1.7)$. We have wall boundary conditions on the bottom, on the top  and on the oblique edge of the mesh, inflow on the left edge and outflow on the right one. The initial conditions are a shock, which divides high density (left--side $x\leq 0$) and low density (right--side $x<0$). This shock has an initial speed in right direction. As the time passes, the shock crosses the oblique surface and creates more internal shock surfaces. The initial conditions are more precisely the following
$$
\begin{pmatrix}
\rho_0\\ u_0\\ v_0 \\ p_0
\end{pmatrix} =
\begin{pmatrix}
8\\
8.25 \\ 0 \\ 116.5
\end{pmatrix} \text{ if } x\leq 0, \qquad \begin{pmatrix}
\rho_0\\ u_0\\ v_0 \\ p_0
\end{pmatrix} =
\begin{pmatrix}
1.4\\
0 \\ 0 \\ 1
\end{pmatrix} \text{ if } x>0.
$$

The parameters used for this test are $\varepsilon=10^{-9}$, convection coefficient $\lambda=15$, CFL = $0.1$, final time $T=0.2$. The mesh we used is composed of $N=19248$ triangular elements with a maximum diameter of $0.0369$.  For $\B^1\, \theta_1=0.1$, for $\B^2\, \theta_1=0.01,\,\theta_2=0.0001$, for $\B^3\, \theta_1=0.005,\,\theta_2=0.0001$.

\begin{figure}[ht!]
	\begin{center}
\includegraphics[width=0.99\textwidth,trim={0 95 0 250},clip]{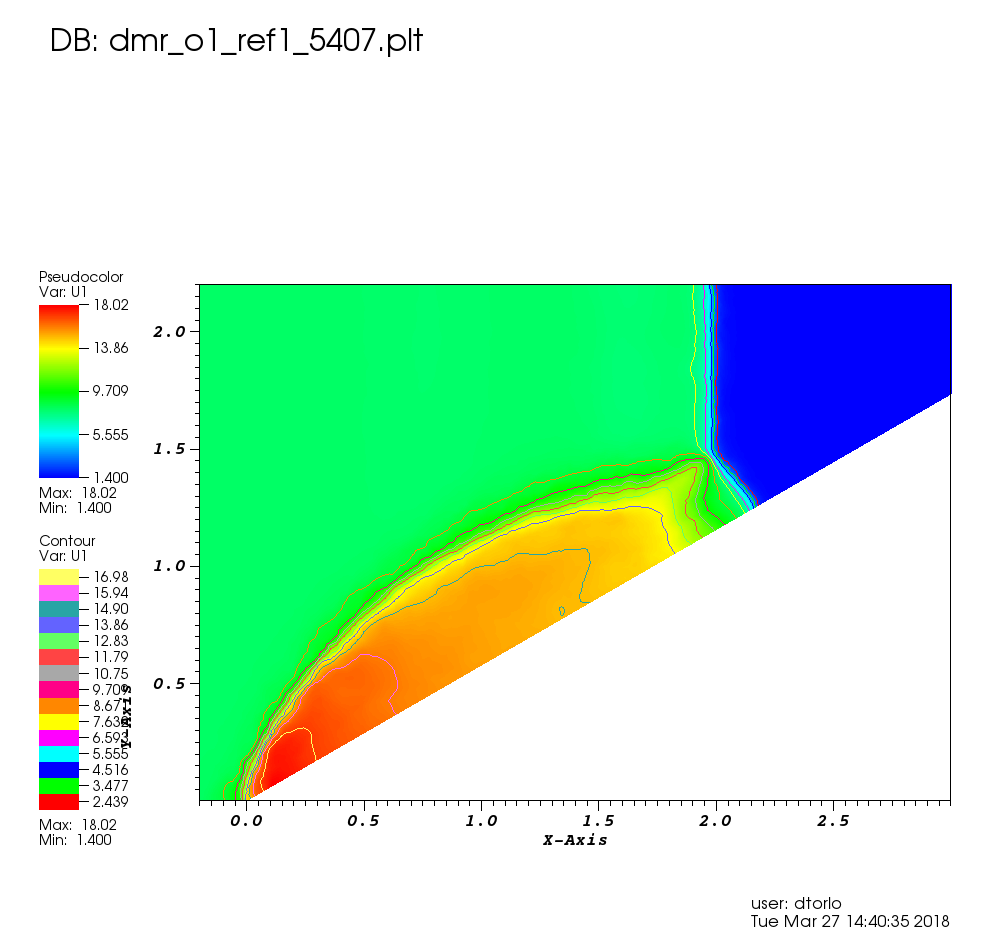}
	\end{center} \caption{Density of DMR test $\mathbb{B}^1$}\label{pic:DMR_2D_surf:B1}
\end{figure}
\begin{figure}[ht!]
	\begin{center}	
\includegraphics[width=0.99\textwidth,trim={0 95 0 250},clip]{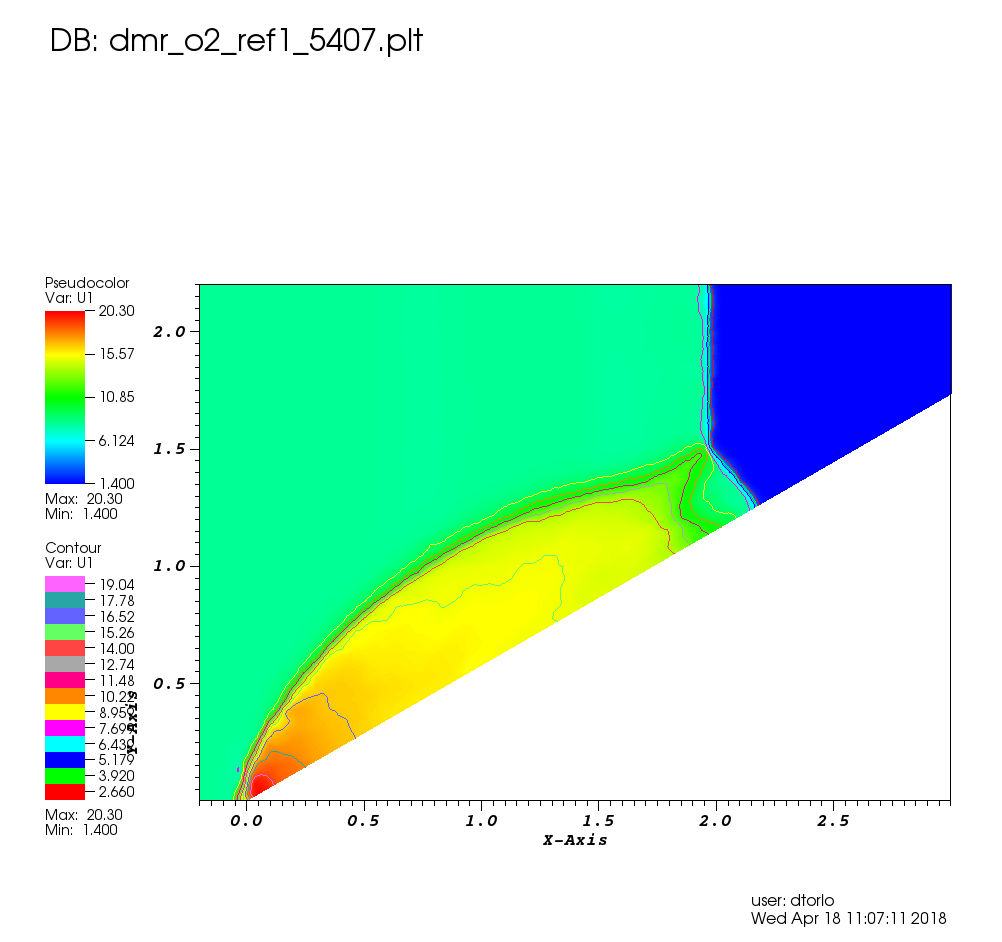}
	\end{center}\caption{Density of DMR test $\mathbb{B}^2$}\label{pic:DMR_2D_surf:B2}
\end{figure}
\begin{figure}[ht!]
	\begin{center}	
\includegraphics[width=0.99\textwidth,trim={0 90 0 250},clip]{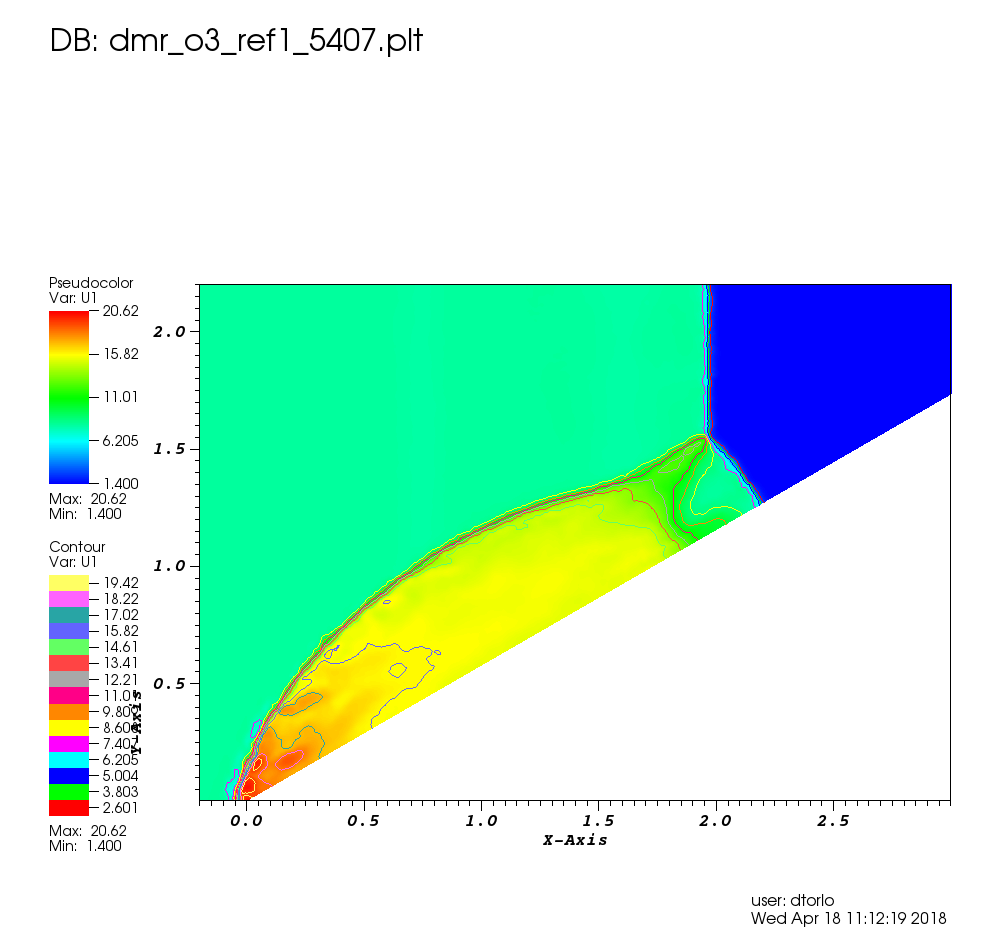}
	\end{center}\caption{Density of DMR test $\mathbb{B}^3$}\label{pic:DMR_2D_surf:B3}
\end{figure}
Again we can see in pictures \ref{pic:DMR_2D_surf:B1}, \ref{pic:DMR_2D_surf:B2} and \ref{pic:DMR_2D_surf:B3} that the scheme catches the behaviour of the shock and its reflection against the lower wall. Even now, we can see that the sharpness of the shock is really well captured by the $\mathbb{B}^3$ scheme, while the others are less precise in defining the shock zone.

\section{Conclusions and further investigations}\label{sec:conclusion}
We have presented a residual distribution scheme for hyperbolic system of equations with stiff relaxation source terms for kinetic models. The method proposed takes advantage of the IMEX formulation (implicit for source term and explicit for advection term) to resolve the stiffness of the relaxation source. Nevertheless, we were able to solve computationally explicitly the kinetic model of \cite{natalini}, thanks to an auxiliary equation, which allows us not to recur to nonlinear solver. The high accuracy of the scheme is reached thanks to two ingredients. 
The first one is the residual distribution framework for spatial discretization \cite{abgrall}, which is a finite element based method that is naturally high order because of the choice of different basis functions. The second is the high order time--integration performed in the DeC method, which allows to couple two schemes, the first easy to solve, for us the IMEX scheme, and a second high order scheme, the high order time--integration residual distribution scheme. The result is an iterative method able to reach high order and stability via few iterations. This is the first time, as far as we know, that the residual distribution framework is used to solve hyperbolic systems with stiff source terms. Even if in this work we solved only one model, it is easily extensible to different models which present similar properties.

The results obtained both from a theoretical point of view and from the simulation side are satisfactory. Indeed, the theorems proved the asymptotic preserving property for our scheme and the rate of accuracy. In addition, the run simulations are reaching the expected accuracy in 1D and 2D, the correct behaviour of the discontinuities of the solutions is well caught by the scheme and as the order increases we can see big improvements in shapes of solutions.

Further investigations may be in the following directions. There are still some open questions over the complete automation of the scheme. For example, it is still not well known which is the relation between parameters $\theta_1, \theta_2$, $\Delta t$ and the quality of the solution. There are studies for 1D smooth solutions, where some relations between these quantities are shown, thanks to some von Neumann stability analysis \cite{stability_RD}. Nonetheless, these results are not easily extensible to nonlinear flux problems or 2D problems.

Moreover, it is not clear why for $\B^3$ the scheme needs more corrections than expected to reach the order of convergence, in particular when the mesh is more refined. This is a contradiction of proposition \eqref{DeC_prop} as shown in \cite{Abgrall2017jcp}.

Finally, we are already working on some extensions of the scheme for multiphase flows equations and we believe that it can be applied also for a large variety of other problems, such as BGK equations, viscoelasticity problems or other kinetic schemes.
\section*{Acknowledgments}
We acknowledge the support of ITN ModCompShock project funded by the European Union’s Horizon 2020 research and innovation program under the Marie Sklodowska-Curie grant agreement No 642768.
We acknowledge Paola Bacigaluppi and Svetlana Tokareva for their contributions in coding and discussing the residual distribution formulation.

\appendix
\section{Residual Distribution schemes}\label{app:residual_distribution}
The key point of the RD schemes is the definition of the splitting of the total residuals into nodal residuals. Through this definition one can actually define the proper scheme to utilise. One can rewrite, for example, the SUPG scheme \cite{hughes} in this way:
\begin{equation}
\phi^K_\sigma (U_h)=\int_{K} \varphi_\sigma ( \nabla A(U_h) - S(U_h)) dx + h_K\int_K \left(\nabla  A(U_h) \cdot \nabla  \varphi_\sigma \right) \tau \left( \nabla  A(U_h) \cdot \nabla  U_h \right).
\end{equation}

What we use in our code are two types of residuals: one for smooth test cases, one for shock test cases. 

\subsection{Smooth solutions residuals}
When we are dealing with smooth tests and we know a priori that we do not need the extra diffusion to dump oscillations brought by discontinuities, we can use a pure Galerkin discretization with a stabilization of jump of the gradient of the solution \cite{burman, DeC_Abgrall}. The study of the stability of the scheme in this situation for smooth solutions is shown in \cite{stability_RD}, through a von Neumann analysis of the scheme.

For an hyperbolic system of equation with source term
\begin{equation}\label{eq:limit_eq_in_append}
\partial_t U + \nabla \cdot A(U) -S(U) = 0,
\end{equation}
the scheme proceeds as follows $\forall \sigma \in \Sigma$
\begin{equation}
\phi_{\sigma}^{K,1} (U_h) = \int_{\partial K} \varphi_\sigma A(U_h)\cdot \textbf{n} d\Gamma -\int_K \nabla \varphi_\sigma \cdot A(U_h) d\textbf{x} -\int_K \varphi_\sigma S(U_h) d\textbf{x},
\end{equation}
and then
\begin{equation}\label{eq:scheme4}
\phi^K_\sigma= \phi^{K,1}_\sigma + \sum_{k=1}^d \sum_{e| \text{edge of }K} \theta_k h_e^{2k} \int_e [\nabla^k U_h ] \cdot [\nabla^k \varphi_\sigma ] d \Gamma.
\end{equation}
Here $d$ is the degree of the polynomial of the basis functions we use, $\theta_k$ are positive coefficients and $[\cdot ]$ is the jump across the edge $e$, namely, if $e$ separates $K$ and $K^+$, $[u]=u|_K-u|_{K^+}$. All the derivatives are meant in the direction of the normal to the edge $e$ and $h_e$ is the length of a 1D element of the mesh (the edge $e$ in 2D, the size of a cell $|K|$ in 1D). The schemes just presented are naturally of order $d+1$ where $d$ is the degree of the polynomial that we are using for the discretization. The parameters $\theta_k$ must be chosen carefully if we want the scheme to be stable. The stability analysis of this scheme in \cite{stability_RD} suggests some optimal values for these parameters in case of 1D linear fluxes. It is not easy to extend this study to different test cases. In addition, these schemes are not too dissipative and they preserve the order of accuracy. Anyway, they do not guarantee stability in case of shocks and discontinuities
 .

\subsection{Shock solutions residuals}
Now, we present the schemes that is used in our simulations in presence of discontinuities or not smooth solutions. More details of these schemes are shown in \cite{paola_svetlana}. The procedure starts defining a local Galerkin Lax--Friedrichs type nodal residual on the steady part of original equation \eqref{eq:limit_eq_in_append}:
\begin{equation}
\phi_{\sigma}^{K,LxF} (U_h): = \int_{\partial K} \varphi_\sigma A(U_h)\cdot \textbf{n} d\Gamma -\int_K \nabla \varphi_\sigma \cdot A(U_h) dx -\int_K \varphi_\sigma S(U_h) dx + \alpha_K (U_\sigma-\overline{U}^K_h),
\end{equation}
where $\overline{U}^K_h$ is the average of $U_h$ over the cell $K$ and $\alpha_K$ is defined as 
\begin{equation}
\alpha_K:=\max\limits_{\sigma\in K} \left( \rho_S  \left( \nabla A(\overline{U}^K_h) \cdot \nabla \varphi_\sigma \right)\right),
\end{equation}
and $\rho_S$ is the function returning the spectral radius of the input matrix. Then, to guarantee monotonicity of the solution near strong discontinuities, we proceed as follows:
\begin{equation}
\begin{split}
&\beta_\sigma ^K (U_h): =\max \left( \frac{\Phi_\sigma^{K,LxF}}{\Phi^K},0 \right) \left( \sum_{j\in K}\max \left( \frac{\Phi_j^{K,LxF}}{\Phi^K},0 \right)\right)^{-1},\\
&\phi^{*,K}_\sigma: =\beta^K_\sigma \phi^K.
\end{split}
\end{equation}
These divisions between vectors are meant component--wise in characteristic variables, that implies the computation of the right eigenvectors of the multiplication of the jacobian of the flux and the normal average velocity $\nabla A(U_h)\cdot \textbf{n}$.
Then, we do a blending between this new residual and the Lax--Friedrichs's one. We use a coefficient $\Theta$ defined as 
\begin{equation}
\Theta:=\frac{|\Phi^K|}{\sum_{j\in K} |\Phi^{K,LxF}_j|}
\end{equation}
and the new residual is 
\begin{equation}
\phi^{\cdot,K}_\sigma:=(1-\Theta)\phi^{*,K}_\sigma + \Theta \Phi^{K,LxF}_\sigma.
\end{equation}
This scheme guarantees the monotonicity principle \cite{abgrall}.

After that we add to the scheme the jump stabilization terms
\begin{equation}\label{eq:scheme5}
\phi^K_\sigma:= \phi^{\cdot,K}_\sigma +\sum_{k=1}^d \sum_{e| \text{edge of }K} \theta_k h_e^{2k} \int_e [\nabla^k U_h ] \cdot [\nabla^k \varphi_\sigma ] d \Gamma,
\end{equation}
and this defines the final scheme.

\section{Deferred Correction properties}
\subsection{Proof of DeC theorem}\label{DeC_proposition}
\begin{proposition}\label{DeC_prop_app}
Let $\L^1$ and $\L^2$ be two operators defined on $\mathbb{R}^m$, which depend on the parameter $\Delta$, such that
\begin{itemize}
\item $\L^1$ is coercive for one norm, i.e., $\exists \alpha_1 >0$ independent of $\Delta$, can be both $\Delta x$ or $\Delta t$ since they are linked by CFL conditions, such that for any $f,g$ we have that $$\alpha_1||f-g||\leq ||\L^1 (f)-\L^1 (g)||$$
\item $\L^1 - \L^2$ is Lipschitz with constant $\alpha_2>0$ uniformly with respect to $\Delta$, i.e., for any $U,V$
$$
||(\L^1(f)-\L^2(f))-(\L^1(g)-\L^2(g))||\leq \alpha_2 \Delta ||f-g||.
$$
\end{itemize}
We also assume that $\exists !\, f^*_\Delta$ such that $\L^2(f^*_\Delta)=0$. Then, if $\eta=\frac{\alpha_2}{\alpha_1}\Delta<1$, the deferred correction is converging to $f^*$ and after $k$ iterations the error is smaller than $\eta^k$.
 
\begin{proof}
Let $f^*$ be the solution of $\L^2(f^*)=0$. Here, we drop the dependency on $f^n$ in $\L^1,\,\L^2$, for simplicity. We know that $\L^1(f^*)=\L^1(f^*)-\L^2(f^*)$, so that
\begin{align}
\L^1(f^{(k+1)})-\L^1(f^*)=&\left(\L^1(f^{(k)})-\L^1(f^*)\right)-\left(\L^2(f^{(k)})-\L^2(f^*)\right),
\end{align} 
then
\begin{align}
\alpha_1 ||f^{(k+1)}-f^*||\leq & ||\L^1(f^{(k+1)})-\L^1(f^*)||=\\
=&||\L^1(f^{(k)})-\L^2(f^{(k)})-(\L^1(f^*)-\L^2(f^*))||\leq \\
\leq & \alpha_2 \Delta ||f^{(k)}-f^*||.
\end{align}
Hence, we can write
\begin{equation}
||f^{(k+1)}-f^*||\leq \left(\frac{\alpha_2}{\alpha_1}\Delta\right) ||f^{(k)}-f^*|| \leq \left(\frac{\alpha_2}{\alpha_1}\Delta\right)^{k+1} ||f^{(0)}-f^*||.
\end{equation}

After $k$ iterations we have an error at most of $\eta^k\cdot ||f^{(0)}-f^*||$. 
\end{proof}
\end{proposition}
\subsection{Lipschitz continuity and coercivity}\label{Dec_properties}
Let us prove that our $\L^1$ and $\L^2$ schemes verify all the hypothesis of proposition \eqref{DeC_prop_app}.
\begin{proposition} $\L^1$ is coercive, i.e., $\exists \alpha_1>0$ s.t. $\forall f,g\in V_h$ and $m =1,\dots, M$
\begin{align}||\L^{1,m}_u(f^0,Pf)- \L^{1,m}_u(f^0,Pg)||&\geq \alpha_1 ||Pf-Pg||,\\
||\L^{1,m}(f^0,f)- \L^{1,m}(f^0,g)||&\geq \alpha_1 ||f-g||.
\end{align}
\begin{proof}
The $u$ part is trivial because 
\begin{equation}
\L^{1,m}_{\sigma , u }(f^0,Pf)- \L^{1,m}_{\sigma , u }(f^0,Pg)= Pf^m_\sigma-Pg^m_\sigma.
\end{equation}
For $f$ part, we have to collect the implicit terms as done in \eqref{eq:L1_f_operator_explicit}. Then, we can write
\begin{equation}
\L^{1,m}_{\sigma }(f^0,f)- \L^{1,m}_{\sigma }(f^0,g)=( f^m_\sigma-g^m_\sigma)-\frac{\Delta t}{\Delta t +\varepsilon}(M(Pf^m_\sigma)-M(Pg^m_\sigma))= f^m_\sigma-g^m_\sigma.
\end{equation}
The last step is possible, since the Maxwellians in our scheme are computed from the $u$ equation and they are actually explicit, so they must coincide. 
If we write the operator explicitly both for $u$ and $f$, we can see that the coercivity constant $\alpha_1=1$, given any norm.
\end{proof}
\end{proposition}

Before proving the Lipschitz continuity, we have to introduce some norms.
We use the following definition of norm for a function $f \in V_h$, which is consistent with the $\L^2$ norm, 
\begin{equation}\label{eq:norm}
||f||^2=\sum_{\sigma \in D_h} |\mathcal{C}_\sigma |f_\sigma^2.
\end{equation}
We also define the norm of all the subtimesteps as
\begin{equation}
|||\mathbf{f}|||=|||(f^0,\dots, f^M)|||=\sqrt{\sum_{m=1}^M ||f^m||^2}.
\end{equation}

Moreover, we will need the definition of the following seminorms
\begin{align}
|f|_{1,x}^2:&= \sum_{\sigma \in D_h}|\mathcal{C}_\sigma |\left( \max_{K | \sigma \in K} \max_{x\in K} \frac{f_\sigma -f(x)}{d(K)}\right)^2,\label{eq:seminorm:space}\\
|\mathbf{f}|_{1,t}^2:&=\sum_{\sigma \in D_h}|\mathcal{C}_\sigma |\left(  \max_{m=1,\dots, M } \frac{f^m-f^{m-1}}{\Delta t^m}\right)^2,\label{eq:seminorm:time}
\end{align}
where $d(K)$ is the diameter of the cell $K$ and it is bounded by $\max_K d(K)=h$. In particular, we note that $|f|_{1,x} \leq |f|_{1}=||\nabla f||_{L^2}$ for every discretization mesh.

\begin{proposition}
If we assume that \begin{align}
&|f|_{1,x} \leq C_1 ||f||,\label{eq:inequality:space_der}\\
&|\mathbf{f}|_{1,t} \leq C_2 |||\mathbf{f}|||,\label{eq:inequality:time_der}
\end{align}
where $C_1 $ and $C_2$ do not depend on the mesh size $h$ and timestep $\Delta t$. And if we require that nodal residuals verify \begin{equation}\label{eq:lipschitz:residuals}
\sum_{\sigma \in D_h}\frac{1}{|\mathcal{C}_\sigma |} \left( \sum_{K|\sigma \in K} \phi^K_\sigma(f) - \phi^K_\sigma (g)\right) ^2 \leq C_3  \sum_{\sigma \in D_h} |\mathcal{C}_\sigma | (f_\sigma-g_\sigma)^2=C_3||f-g||^2.
\end{equation}
Then, $\L^1-\L^2$ is Lipschitz continuous, i.e., $\exists \alpha_2>0$  s.t. $\forall f,g \in V_h$ 
\begin{align}
\label{eq:lipschitz_u}
||| \left(\L^{1}_u(P\mathbf{f})- \L^{1}_u(P \mathbf{g}) \right)- \left(\L^{2}_u(P\mathbf{f})- \L^{2}_u(P\mathbf{g}) \right)||| &\leq \alpha_2 \Delta |||P\mathbf{f}-P\mathbf{g}|||,\\ 
||| \left(\L^{1}(\mathbf{f})- \L^{1}(\mathbf{g})\right) - \left(\L^{2}(\mathbf{f})- \L^{2}(\mathbf{g}) \right)||| &\leq \alpha_2 \Delta |||\mathbf{f}-\mathbf{g}|||. \label{eq:lipschitz_f}
\end{align}

\paragraph{Remark}
The extra hypothesis added are related to the regularity of the solution. Of course, this is not always the case, and, for example, when there are shocks in the solution, \eqref{eq:inequality:space_der} does not hold. Anyway, even if we can not prove the convergence for those cases, we see numerically a big improvement in higher order solutions.
The inequality \eqref{eq:inequality:time_der} is actually given, during the DeC procedure, by the Lipschitz continuity of fluxes and residuals. To keep the proof more general, we add it as a further hypothesis. Equation \eqref{eq:lipschitz:residuals}, in our case, is given by the consistency of the nodal residuals, the Lipschitz continuity of the flux $F$ and by the regularity of the solutions $f,g$ as stated in \eqref{eq:inequality:space_der}.

\begin{proof}
The estimation of \eqref{eq:lipschitz_u} is a simplification of the case of \eqref{eq:lipschitz_f}, so we will skip its proof. 

For simplicity, let us define the differences $\delta f:=f-g$, $\delta \phi^K_\sigma(f):=\phi^K_\sigma (f)-\phi^K_\sigma (g)$, $\delta M(Pf):= M(Pf)-M(Pg)$, $\delta \L:= \L^1-\L^2$ and $\delta \mathcal{I}(\mathbf{f}):=\mathcal{I}_0(\mathbf{f})-\mathcal{I}_M(\mathbf{f})$.

Let us split the operators into two parts. The first one is composed of the term related to time derivative and source term $\L_{ts}$, the second one concerns the advection part $\L_{ad}$. If we write explicitly the source and time part, we get
\begin{subequations}
\begin{equation}
\begin{split}
\delta \L^{m}_{ts,\sigma} (f)-&\delta \L^{m}_{ts,\sigma} (g)   =\\
= \sum_{K|\sigma \in K} \frac{1}{|\mathcal{C}_\sigma |}\bigg[ &\frac{\varepsilon}{\varepsilon+\Delta t^m}  \int_K \varphi_\sigma \left(\delta f^m_\sigma  - \delta f^m \right) - \frac{\Delta t}{\varepsilon}\int_K \varphi_\sigma \left( \delta M(Pf^m_\sigma) - \delta f^m_\sigma \right)+\\
+& \frac{\varepsilon}{\varepsilon+\Delta t^m} 
\frac{1}{\varepsilon} \int_{t^0}^{t^m} \mathcal{I}_M \left( \delta \phi^K_{s,\sigma} ( f^0), \dots ,  \delta \phi^K_{s,\sigma} ( f^M) ,s \right) ds \bigg] .
\end{split}
\end{equation}
\end{subequations}
Now, let us suppose that the residuals are a consistent discretization of fluxes and source terms, so let us use the Galerkin discretization instead of any other one. Moreover, let us add and subtract the residual in timestep $t^{n,m}$. So, we can write, neglecting $\mathcal{O}(\Delta ^2)$,
\begin{subequations}
\begin{align}
& \L^{1,m}_{ts,\sigma} (f)-\L^{1,m}_{ts,\sigma} (g)- \L^{2,m}_{ts,\sigma} (f)+\L^{2,m}_{ts,\sigma} (g) +\mathcal{O}(\Delta ^2) = \\
\begin{split}
=&\frac{1}{|\mathcal{C}_\sigma |}\int_\Omega \varphi_\sigma \left(  \delta f^m_\sigma - \delta f^m   \right)-  \frac{1}{|\mathcal{C}_\sigma |}\frac{\Delta t^m}{(\varepsilon+\Delta t^m)}  \int_\Omega \varphi_\sigma \left(\delta M(Pf^m_\sigma)-   \delta M(Pf^m)  \right) +\\
+& \frac{1}{\varepsilon+\Delta t^m} \int_{t^{0}}^{t^m} \mathcal{I}_M(\delta \phi^K_{s,\sigma}(f^0)-\delta \phi^K_{s,\sigma}(f^m), \dots,\delta \phi^K_{s,\sigma}(f^M)-\delta \phi^K_{s,\sigma}(f^m),s)ds .
\end{split}\label{eq:f_collecting}
\end{align}
\end{subequations}
\begin{subequations}
Now, we sum over the DoFs and we square the previous quantity. We use Lemma A.1 of \cite{DeC_Abgrall} to pass from coefficients $v_\sigma$ to pointwise evaluation $v(\sigma)$, with the abuse of notation. It states that $\sum_{\sigma\in K} |v_\sigma - v_{\sigma '}| \leq C_K \sum_{\sigma \in K}|v(\sigma)-v(\sigma')|$ where $C_K$ is the norm of the inverse of the matrix $(\varphi_\sigma(\sigma'))_{\sigma, \sigma'}$ and it depends on $K$ only via the aspect ratio of the element $K$.  
\begin{align}
&\sum_{\sigma\in D_h}|\mathcal{C}_\sigma | \left( \L^{1,m}_{ts,\sigma} (f)-\L^{1,m}_{ts,\sigma} (g)- \L^{2,m}_{ts,\sigma} (f)+\L^{2,m}_{ts,\sigma} (g) \right )^2  \leq \\
\begin{split}\label{eq:scale_outside}
\leq &C_a h^2 \sum_{\sigma \in D_h }  \frac{1}{|\mathcal{C}_\sigma |}\left( \int_\Omega \varphi_\sigma \left( \frac{ \delta f^m_\sigma  - \delta f^m(x)}{d(K)}  \right) \right)^2 + \\ 
+&C_b h^2  \frac{\Delta t ^m}{ (\varepsilon+\Delta t^m)}\sum_{\sigma \in D_h }  \frac{1}{|\mathcal{C}_\sigma |}\left(\int_\Omega \varphi_\sigma \frac{ \delta M(Pf^m)(\sigma) -\delta M(Pf^m)}{d(K)} \right)  ^2+\\
+&C_c\frac{\Delta t^m}{\varepsilon + \Delta t^m}\sum_{\sigma\in D_h} |C_\sigma| \max_{r} \left( \delta \phi^K_{s,\sigma}(f^r)-\delta\phi^K_{s,\sigma}(f^m) \right)^2  \leq
\end{split}\\ \label{eq:gradient}
\begin{split}
\leq &C_d h^2  (|\delta f^m|_{1,x}^2 + |\delta M(Pf^m)|_{1,x}^2+ \max_{r} ||\delta f^r - \delta f^m|| ) \leq
\end{split}\\
\label{eq:seminorms}
\leq &C_e h^2  (||\delta f^m||^2 + ||\delta M(Pf^m)||^2 + \Delta t^2 |\delta \mathbf{f}|_{1,t}) \leq \\
 \label{eq:lipschitz_maxwellians}
\leq &C_f h^2  |||\delta \mathbf{f}|||^2  + \mathcal{O}(h^4) \leq C_4 h^2 |||\mathbf{f-g}|||^2.
\end{align}
\end{subequations}
In \eqref{eq:scale_outside} we explicitly bring the scale $h$ outside the first two sums, while in the third term we just bound the interpolant polynomail with the maximum of the interpolant values times a constant, in \eqref{eq:gradient} we use the definition of the seminorm \eqref{eq:seminorm:space}, the Lipschitz continuity of residuals \eqref{eq:lipschitz:residuals}, the product rule for integrals and the bound $\Delta t^m\leq \Delta t^m +\varepsilon$. 
In  \eqref{eq:seminorms} we use the inequality \eqref{eq:inequality:space_der} and the definition of the seminorm \eqref{eq:seminorm:time}. In \eqref{eq:lipschitz_maxwellians} we use the fact that the maxwellians $M$ and the projections $P$ are Lipschitz continuous, the inequality \eqref{eq:inequality:time_der} and the fact that $\Delta t \sim h$. 
The constant $C_4$ does not depend on $h, \Delta t$ nor on $\varepsilon$. It depends on the size of the domain, on the Lipschitz continuity of the Maxwellians, on the regularity of the mesh and on basis functions.

For the advection term a similar computation is carried out, but, in this case the error is a $\mathcal{O}(\Delta t)$. Using the notation of $\phi_\sigma:=\sum\limits_{K|\sigma \in K} \phi^K_\sigma$, let us write 
\begin{subequations}
\begin{align}
&||\mathcal{S}_x||^2:=\sum_{\sigma \in D_h}  |\mathcal{C}_\sigma |\left( \delta \L^{1,m}_{ad,\sigma} (f)-\delta \L^{1,m}_{ad,\sigma} (g) \right) ^2  = \\
\begin{split}
=&\sum_{\sigma \in D_h}  \frac{1}{|\mathcal{C}_\sigma |} \bigg( \frac{\varepsilon}{\varepsilon + \Delta t^m} \int_{t^{n,0}}^{t^{n,m}} \delta \mathcal{I} \left( \delta\phi_{ad,\sigma} ( f^0) , \dots , \delta \phi_{ad,\sigma} ( f^M),s \right) ds\bigg)^2 \leq 
\end{split}\label{eq:f_advection} \\
\leq &  C_l \sum_{\sigma \in D_h}   \frac{\Delta t^2}{|\mathcal{C}_\sigma |} \left(\sum_{K|\sigma \in K} \max_{m=1,\dots, M} \frac{ |\delta \phi^K_{ad,\sigma} (f^m)-\delta \phi^K_{ad, \sigma} (f^{m-1})|}{\Delta t^m}  \right)^2.  \label{eq:adv_interpolant} 
\end{align}
In \eqref{eq:adv_interpolant} we use the bound $\varepsilon\leq \varepsilon + \Delta t^m$ and the fact that $\mathcal{I}_0$ is a zero order approximation of $\mathcal{I}_M$, so, adding the integration in time, we get the error estimation above. 
\begin{align}
||\mathcal{S}_x||^2 \leq C_q \sum_{\sigma \in D_h}  & \Delta t^2|\mathcal{C}_\sigma | \left(\max_{m=1,\dots, M} \frac{ |\delta f^m-\delta f^{m-1}|}{\Delta t^m}  \right)^2\leq \label{eq:f:lipschitz_residual}\\
\leq C_p \Delta t^2   &  \sum_{m=1}^{M}| f^m-g^m|_{1,t}^2 \leq C_5 \Delta t^2|||\mathbf{f}-\mathbf{g}|||^2. \label{eq:f:bound:seminorm:time}
\end{align}

In \eqref{eq:f:lipschitz_residual} we use the Lipschitz continuity and consistent hypothesis over the residuals as stated in \eqref{eq:lipschitz:residuals}. Finally, in \eqref{eq:f:bound:seminorm:time} we use the definition of seminorm \eqref{eq:seminorm:time} and we apply the bound in \eqref{eq:inequality:time_der}. Again, $C_5$ does not depend on $\Delta t$, $h$ or $\varepsilon$, but only on fluxes, geometry and basis functions.

In conclusion, summing up the inequalities \eqref{eq:lipschitz_maxwellians} and \eqref{eq:f:bound:seminorm:time}, we prove the thesis of the proposition.
\end{subequations}
\end{proof}
\end{proposition}

\clearpage
\bibliographystyle{plain}

\bibliography{../../../biblio}

\end{document}